# A UNIFIED FRAMEWORK FOR UTILITY MAXIMIZATION PROBLEMS: AN ORLICZ SPACE APPROACH

By Sara Biagini and Marco Frittelli

*Università degli Studi di Perugia and Università degli Studi di Milano*

We consider a stochastic financial incomplete market where the price processes are described by a vector-valued semimartingale that is possibly nonlocally bounded. We face the classical problem of utility maximization from terminal wealth, with utility functions that are finite-valued over $(a, \infty)$, $a \in [-\infty, \infty)$, and satisfy weak regularity assumptions. We adopt a class of trading strategies that allows for stochastic integrals that are not necessarily bounded from below. The embedding of the utility maximization problem in Orlicz spaces permits us to formulate the problem in a unified way for both the cases $a \in \mathbb{R}$ and $a = -\infty$. By duality methods, we prove the existence of solutions to the primal and dual problems and show that a singular component in the pricing functionals may also occur with utility functions finite on the entire real line.

**1. Introduction.** In the most general semimartingale model for the underlying process $S$, the problem we address takes the form

$$\sup_{H \in \mathcal{H}} E[u(x + (H \cdot S)_T)], \tag{1}$$

where:

- $u : \mathbb{R} \to \mathbb{R} \cup \{-\infty\}$ is the utility function of the agent, which is assumed to be increasing and concave on the interior $(a, \infty)$, $a \in [-\infty, \infty)$, of its effective domain and to satisfy $\lim_{x \to -\infty} u(x) = -\infty$;
- $x > a$ is the initial endowment of the agent and $T \in (0, \infty]$ is the time horizon;
- the process $S$ is an $\mathbb{R}^d$-valued càdlàg semimartingale defined on the filtered probability space $(\Omega, \mathcal{F}, (\mathcal{F}_t)_{t \in [0,T]}, P)$, the filtration satisfies the usual assumptions of right continuity and completeness and $\mathcal{F}_0$ is *trivial*, that is, it is generated by the $P$-negligible sets in $\mathcal{F}_T$; in case $T = \infty$, we assume that, for every process $Y$ considered, the limit $Y_\infty = \lim_{t \uparrow +\infty} Y_t$ exists;
- $\mathcal{H}$ is a class of admissible $\mathbb{R}^d$-valued $S$-integrable predictable processes, which represents the allowed trading strategies [see (4) for the precise definition];
- $H \cdot S$ is the stochastic integral and $(H \cdot S)_T$ is the terminal gain achieved by following the strategy $H$.

The "duality approach" to the resolution of this very classical problem was first employed by [20] (see also [10] for earlier work in stochastic optimal control) and is based on classical tools from convex analysis. As far as we know, we give the most general formulation of the duality to date.

To formulate the dual optimization problem, we denote by $\Phi : \mathbb{R}_+ \to \mathbb{R} \cup \{\infty\}$ the function

$$\Phi(y) \triangleq \sup_{x \in \mathbb{R}} \{u(x) - xy\},$$

which is the convex conjugate of the utility function $u$. The dual problem for the utility maximization is typically

$$(2) \qquad \min_{\lambda > 0, Q \in \mathcal{M}} \lambda x + E\left[\Phi\left(\frac{dQ}{dP}\right)\right],$$

where $\mathcal{M}$ is an appropriate set of measures, but, under our assumptions, it will be a generalized form of (2). In Section 3, we will see that the dual variables are not only probabilities, but possibly more general functionals.

We set

$$\mathbb{P}_\Phi \triangleq \left\{ Q \ll P \,\Big|\, E\left[\Phi\left(\frac{dQ}{dP}\right)\right] < \infty \right\}.$$

The following assumptions will not be needed until Section 4, but it is worthwhile to formulate them here so that appropriate comparison with existing literature is possible.

(A1) The utility function $u : \mathbb{R} \to \mathbb{R} \cup \{-\infty\}$ is increasing, strictly concave and continuously differentiable on the interior $(a, \infty)$, $a \in [-\infty, \infty)$, of its effective domain and satisfies the Inada conditions

$$u'(a) \triangleq \lim_{x \downarrow a} u'(x) = +\infty, \qquad u'(\infty) \triangleq \lim_{x \uparrow \infty} u'(x) = 0.$$

(A2)

$$(3) \qquad \mathbb{P}_\Phi = \mathbb{P}_{\Phi_\lambda},$$

where the function $\Phi_\lambda : \mathbb{R}_+ \to \mathbb{R}$ is defined by $\Phi_\lambda(y) \triangleq \Phi(\lambda y)$, with $\lambda > 0$ fixed.

REMARK 1. The condition (3) involves not only the function $u$ (through its conjugate function $\Phi$) but also the probability measure $P$. When the probability space is finite and $\Phi$ is finite-valued on $(0, \infty)$, (3) is always satisfied, regardless of the growth properties of $\Phi$. In [5], Section 2.2, we showed that (3) is weaker than the condition of reasonable asymptotic elasticity on $u$ [RAE($u$)] introduced by Schachermayer [24]. On the relationship between RAE($u$), condition (3) and the $\Delta_2$-condition in Orlicz space theory, see [19], Section 6, for the case where $a$ is finite and [24] or [5], Section 2.2, for the case where $a = -\infty$.

We now discuss the literature that considers the utility maximization problem in the context of $\mathbb{R}^d$-valued semimartingale price processes and that is not restricted to a particular utility function. The interested reader may find exhaustive references in [5, 19, 24].

The current literature is essentially split into two main branches.

1. First case: $a \in \mathbb{R}$, so that the utility functions have a half-line as proper domain, for example, $u(x) = \sqrt{x-a}$, $u(x) = \ln(x-a)$.

    Under (A1) and the assumption that the asymptotic elasticity of $u$ at $+\infty$ is strictly smaller than 1 [$AE_{+\infty}(u) < 1$], and when $S$ is a general semimartingale, this subject was thoroughly analyzed in [19] and [12]. In the first paper, the assumption $AE_{+\infty}(u) < 1$ was introduced and it was shown to be crucial for the existence of the solution of problem (1). As shown in Remark 39 of Section 6.1, when $a$ is finite, the condition $AE_{+\infty}(u) < 1$ implies (A2).

    In the cited references, it was also shown that the dual variables $Q \in \mathcal{M}$ may not be true probabilities and a singular component may show up. This is particularly evident in the approach of [12], where $\mathcal{M} \subseteq ba(\Omega, \mathcal{F}, P)$, the space of finitely additive measures on $\mathcal{F}$ that are absolutely continuous with respect to $P$. These authors also remarked that the solution of the dual problem may not be unique, but no explicit example was given.

2. Second case: $a = -\infty$, so that the utility functions have $\mathbb{R}$ as proper domain, for example, $u(x) = -e^{-\gamma x}, \gamma > 0$.

   - Under assumptions (A1) and RAE($u$) [stronger than (A2)], and when $S$ is a *locally bounded semimartingale*, the problem was addressed in [24]. The set $\mathcal{H}$ of strategies employed here is the classical set $\mathcal{H}^1$ of strategies with wealth uniformly bounded from below. The dual problem has exactly the form (2) and the dual variables are local martingale probabilities for $S$.

     As regards the optima, one cannot expect the solution of the primal problem to be bounded from below, so that, in general, it will not

belong to $\mathcal{H}^1$. However, in [24], the problem was reformulated by considering the $L^1(P)$-closure of the set of random variables $u(G)$, where $G \in L^0$ is dominated by some terminal gain $(H \cdot S)_T$, with $H \in \mathcal{H}^1$. Under the assumption RAE($u$), it was shown in [24]—in addition to several key duality results—that the primal solution $f_x$ exists in this enlarged set and that the dual solution $Q_x$ is unique. In addition, it was proven that $f_x$ can be represented as a stochastic integral, as soon as $Q_x$ is equivalent to $P$.

- Under assumptions (A1) and (A2) and when $S$ is a *general*, *possibly nonlocally bounded*, *semimartingale*, the problem was analyzed in [5]. We discuss these results in detail, since we will be adopting some of the definitions and notation.

The work [5] is based on a careful analysis of the proper set of strategies $\mathcal{H}$ that are allowed in the trading. Indeed, the traditional set of strategies $\mathcal{H}^1$ may reduce to the null strategy when $S$ is nonlocally bounded, so the maximization problem on this set may turn out to be trivial. This may happen if, for example, $S$ is a compound Poisson with unbounded jump size (see also the toy Example 4 below).

To model the situation in which the investor is willing to take more risk to really increase his/her expected utility in a very risky market, in [5], we enlarged the set of allowed strategies by admitting losses bounded from below by $-cW$, where $W \geq 0$ is a random variable, possibly unbounded from above. We defined the set $\mathcal{H}^W$ of $W$-admissible strategies by

(4) $\qquad \mathcal{H}^W = \{H \in L(S) \mid (H \cdot S)_t \geq -cW \ \forall t \leq T, \text{ for some } c > 0\},$

where $L(S)$ is the class of predictable and $S$-integrable processes. We showed that the stochastic integrals associated with these strategies enjoy good mathematical properties when the random variable $W$ that controls the losses satisfies the conditions of suitability for the market and compatibility with the preferences. Here are the definitions.

DEFINITION 2. $W \in L^0_+$ is *suitable* (for the process $S$) if $W \geq 1$ and, for each $i = 1, \ldots, d$, there exists $H^i \in L(S^i)$ such that

$$P(\{\omega \mid \exists t \geq 0 \ H^i_t(\omega) = 0\}) = 0$$

and

(5) $\qquad |(H^i \cdot S^i)_t| \leq W \qquad \text{for all } t \in [0, T], P\text{-a.s.}$

The set of suitable random variables is denoted by $\mathbb{S}$.

DEFINITION 3. $W \in L^0_+$ is *compatible* (with the preferences of the agent) if

(6) $\qquad \forall \alpha > 0 \qquad E[u(-\alpha W)] > -\infty.$

Notice that $\mathcal{H}^W = \mathcal{H}^{\alpha W}$ for all $W \geq 0$ and constants $\alpha > 0$, so that $\mathcal{H}^W$ does not change if $W$ is scaled by a multiplicative factor. Therefore, the request $W \geq 1$ in the definition of suitability is only intended to guarantee that $W$ is bounded away from zero.

When $S$ is locally bounded, $W = 1$ is automatically suitable and compatible (see [5], Proposition 1), while, in general, there is no natural choice for $W$ (if there is any, see Example 4).

In Section 3.1, the role played by the *suitability* condition in ensuring that the (regular) dual variables are $\sigma$-martingale measures will become evident.

Under assumptions (A1) and (A2), we then proved the subsequent results.

(a) For all loss variables $W$ that are compatible and suitable, the optimal value on the class $\mathcal{H}^W$ coincides with the optimal value $U_\Phi(x)$ of the maximization problem over a larger domain $K_\Phi$. The set $K_\Phi$ and the value $U_\Phi(x)$ do not depend on the single $W$, but depend on the utility function $u$ through its conjugate function $\Phi$.

(b) For all loss variables $W$ that are compatible and suitable, the following duality relation holds true:

$$\text{(7)} \quad \sup_{H \in \mathcal{H}^W} E[u(x + (H \cdot S)_T)] \\ = \min_{\lambda > 0, Q \in \mathbb{M}_\sigma \cap \mathbb{P}_\Phi} \left\{ \lambda x + E\left[\Phi\left(\lambda \frac{dQ}{dP}\right)\right] \right\} := U_\Phi(x),$$

where $\mathbb{M}_\sigma$ is the set of $\sigma$-martingale measures absolutely continuous with respect to $P$.

(c) The primal solution $f_x$ exists in the set $K_\Phi$, but, in general, it does not belong to $\{(H \cdot S)_T \mid H \in \mathcal{H}^W\}$. The dual solution $Q_x$ is unique.

(d) $f_x$ is $Q_x$-a.s. equal to the terminal value of a stochastic integral $(H_x \cdot S)_T$.

The pleasing property that the dual variables are probabilities is, in fact, ensured by the compatibility condition (6), as will be clarified in Section 6.2.

However, it is clear that (6) puts some restrictions on the jumps of $S$, as highlighted in the next toy example.

EXAMPLE 4. Consider a single period market model with $S_0 = 1$ and trivial initial $\sigma$-algebra $\mathcal{F}_0$. Let $(\Omega, \mathcal{F}_1, P) = (\mathbb{R}, \mathcal{B}(\mathbb{R}), \psi(x)\, dx)$, where $dx$ is the Lebesgue measure and $\psi$ is a density function on $\mathbb{R}$, and let $S_1 : \mathbb{R} \to \mathbb{R}$ be the identity map. Then, $S = (S_0, S_1)$ is a semimartingale which is nonlocally bounded as soon as the support of $\psi$ is unbounded. Let us assume that $S$ is nonlocally bounded and note that, $\mathcal{H}^1 = \{0\}$, so the constant 1 is not suitable. In this model, it is easy to see that, basically, the unique suitable $W$ is $W = 1 + |S_1|$ and, consequently, $\mathcal{H}^W = \mathbb{R}$. Let us select the exponential utility $u(x) = -e^{-x}$ and check the compatibility of $W = 1 + |S_1|$ in the situations below.

1. If $\psi$ is a Gaussian density, then $W$ satisfies the compatibility condition (6).
2. If $\psi$ is a two-sided exponential density [e.g., $\psi(x) = \frac{\lambda}{2} e^{-\lambda |x|}$], then $W$ does not verify (6), since

$$E[u(-\alpha W)] = E[-e^{\alpha W}] > -\infty$$

holds true only if $0 \leq \alpha < \lambda$.
3. If $\psi$ is a Cauchy density [e.g., $\psi(x) = \frac{1}{\pi(1+x^2)}$], then we have an extreme case in which $E[u(-\alpha W)] > -\infty$ only if $\alpha = 0$. The expected utility from nonzero investments in $S$ is always $-\infty$.

 Informally speaking, in this case, the exponential utility is totally incompatible with the market structure.

As the extreme "incompatibility" case in item 3 above shows, a reasonable utility maximization problem cannot be built without *any* restrictions. But, at the same time, this example suggests that condition (6) may be relaxed to cover new, interesting situations similar to that in item 2 of the toy model.

We thus introduce a milder notion of compatibility.

DEFINITION 5. $W \in L_+^0$ is weakly *compatible* (with the preferences of the agent) if

(8) $\qquad \exists \alpha > 0 \qquad E[u(\alpha W)] > -\infty.$

We can finally list the main results of this paper (see Theorems 21 and 29).

- We simultaneously treat the cases $a$ finite and $a = -\infty$.
- We extend the aforementioned results of [5] by adopting condition (8) on $W$, which allows us to consider more general market models.
- We prove that a duality relation holds and show that, in general (*even with exponential utility functions*), the solution of the dual problem will have a singular component, in a sense to be clarified in Section 2.2.
- Under the assumptions (A1) and (A2), we show that the primal solution $f_x$ exists in an enlarged set $K_\Phi^W$ and we characterize it in terms of the dual solution.
- Under the assumption of the existence of a suitable and *compatible* loss variable, we prove (in Section 7) that the optimal value on the class $\mathcal{H}^W$ for *any weakly compatible* (resp. compatible) $W \in L_+^0$ is bounded by (resp. equal to) $U_\Phi(x)$. Hence, under this assumption, there is no incentive to enlarge the set of strategies by adopting a weakly compatible $W$.

As shown in Example 4, the set of suitable and *compatible* loss variables may be empty. In this case, the optimal level of wealth from the class $\mathcal{H}^W$ may depend on the selection of the particular weakly compatible $W$. A thorough study of this issue is left for future research.

The occurrence of the singular component is a consequence of the potentially big losses admitted in trading. In Section 4.1, we prove that under some circumstances—when the optimal loss is well inside the tolerated margin—the singular component is again zero.

As regards the representation of the optimal $f_x$ as terminal value of a stochastic integral process and the supermartingale property of this integral process, we refer to [24] for the locally bounded case and to [5, 6] when $S$ is a general semimartingale and the set of suitable and *compatible* loss bounds is not empty. The extension of these results when the loss bound $W$ can only be *weakly compatible* is left for future investigation.

We think that one major novelty of the paper is the first point of the above list, which is rather a philosophical contribution and thus more valuable.

*We believe that there are no good reasons for treating the problem* (1) *separately—for the two cases $a = -\infty$ or $a$ finite—as has been done until now. These two apparently different situations can, in fact, be seen as particular cases of a single, unified framework.*

In this paper, the definitions of admissible trading strategies and the domains of the primal and dual optimization problems are the same for both cases. Moreover, the proofs of the main results are all formulated in the unified setup.

In Section 6, we also show how the known results in [12] and [5] can be deduced as corollaries of our theorems.

In Section 2, we will introduce the duality framework that is the key tool in our unified presentation.

In [9], it was first shown how to use Orlicz space duality to address the utility maximization problem in the case where $a = -\infty$ and $W$ satisfies condition (6). The key point there is that the Orlicz spaces in question are naturally induced by the utility function $u$.

Following the ideas in [9], in Section 2, we build a duality framework, which also works for $a$ finite and when $W$ satisfies condition (8).

The basic idea behind the construction of the Orlicz duality is the following. Given that the utility function $u$ is concave, the wildest behavior is seen on the left tail, that is, the losses are weighted in a more severe way than the gains. This simply reflects the risk aversion of the agent. The left tail of $u$ can easily be turned into the Young function $\widehat{u}(x) := -u(-|x|) + u(0)$, thus giving rise to an appropriate Orlicz space. Condition (8) then just means that $W$ (which is positive) belongs to the Orlicz space $L^{\widehat{u}}(P)$, while condition (6) would mean that $W$ belongs to a "good" subspace of $L^{\widehat{u}}(P)$. By the definition of $\widehat{u}$, only the negative values of $u$ are taken into account and the Orlicz space $L^{\widehat{u}}(P)$ will only be used *to control the possible losses* occurring in trading.

The rest of the article is organized as follows. After some preliminary results, in Section 3, we introduce the set of dual variables and state the duality theorem (Theorem 21). The existence of the solution to the primal problem and some of its properties are proved in Section 4. In Section 5, we give some examples. In particular, we describe a concrete example in the case $a = -\infty$ where there are infinite dual solutions, explicitly characterized.

**2. The Orlicz spaces associated with $u$ and $\Phi$.** Up to Section 4, the utility functions $u : \mathbb{R} \to \mathbb{R} \cup \{-\infty\}$ are increasing, concave on the interior $(a, \infty)$ of the effective domain and satisfy $\lim_{x \to -\infty} u(x) = -\infty$ [it is understood that $u(x) = -\infty$ for all $x < a$, if $a$ is finite]. Without loss of generality, we may (and do) assume, from now on, that

$$a < 0$$

(this can always be obtained by translation if $a$ is finite). From our results in the case where $a$ is negative and finite, one may easily recover the corresponding results (see Section 6.1) in the case $a \geq 0$.

Under these conditions, the function $\Phi$ conjugate to $u$ satisfies $\Phi(\infty) = \infty$ and $\Phi(0^+) = u(\infty)$.

Two widely used utility functions that satisfy the above requirements are the logarithm $u(x) = \ln(1 + x)$ ($a = -1$) and the exponential $u(x) = -e^{-x}$ ($a = -\infty$). But our class of utility functions includes functions for which $\lim_{x \downarrow -\infty} \frac{u(x)}{x}$ is finite (positive) and/or functions constant for $x \geq x_0$.

In the paper, the $L^p(\Omega, \mathcal{F}, P)$ spaces, $p = 0$ or $p \in [1, \infty]$, will simply be denoted by $L^p$, unless it is necessary to specify the probability, in which case we write $L^p(P)$.

In Section 2.1, we recall the generalities on Orlicz spaces and introduce the appropriate Orlicz spaces $L^{\widehat{u}}, L^{\widehat{\Phi}}$, which are constructed with the same methodology of [9], thanks to our assumption $a < 0$.

Section 2.2 deals with $(L^{\widehat{u}})^*$, the norm dual of $L^{\widehat{u}}$, and its decomposition. We pay particular attention to the *singular elements* of $(L^{\widehat{u}})^*$ and their properties.

2.1. *Generalities.* A Young function $\Psi$ is an even, convex function $\Psi : \mathbb{R} \to \mathbb{R} \cup \{+\infty\}$ with the following properties:

1. $\Psi(0) = 0$;
2. $\Psi(\infty) = +\infty$;
3. $\Psi < +\infty$ in a neighborhood of 0.

Note that $\Psi$ may jump to $+\infty$ outside of a bounded neighborhood of 0. In case $\Psi$ is finite-valued, however, it is also continuous, by convexity.

The Orlicz space $L^\Psi(P)$, or simply $L^\Psi$, on $(\Omega, \mathcal{F}, P)$ is then defined as

$$L^\Psi = \{f \in L^0 \mid \exists \alpha > 0 \ E[\Psi(\alpha f)] < +\infty\}.$$

It is a Banach space with the Luxemburg (or gauge) norm

$$N_\Psi(f) = \inf\left\{c > 0 \ \Big| \ E\left[\Psi\left(\frac{f}{c}\right)\right] \leq 1\right\}.$$

With the usual pointwise lattice operations, $L^\Psi$ is also a Banach lattice, that is, the norm satisfies the monotonicity condition

$$|g| \leq |f| \Rightarrow N_\Psi(g) \leq N_\Psi(f).$$

It is not difficult to prove that

$$L^\infty \hookrightarrow L^\Psi \hookrightarrow L^1$$

with linear lattice embeddings (the inclusions). In fact, these spaces are a generalization of the familiar $L^p$ spaces. To recover $L^p$ with $1 \leq p < +\infty$, take $\Psi_p(x) = |x|^p$ as Young function. To recover $L^\infty$, consider the Young function $\Psi_\infty(x) = \delta_C(x)$, where $\delta_C$ is the indicator function of the convex set $C = \{x \in \mathbb{R} \mid |x| \leq 1\}$ ($\delta_C = 0$ on $C$ and $\delta_C = +\infty$ on $\mathbb{R} \setminus C$).

There is an important linear subspace of $L^\Psi$, namely

$$M^\Psi = \{f \in L^0 \mid E[\Psi(\alpha f)] < +\infty \ \forall \alpha > 0\}.$$

In general, $M^\Psi \subsetneq L^\Psi$. This can be easily seen when $\Psi = \Psi_\infty$ since, in this case, $M^{\Psi_\infty} = \{0\}$, but there are also nontrivial examples of the strict containment with finite-valued, continuous Young functions that we will consider soon.

However (see [23]), when $\Psi$ satisfies the following $\Delta_2$ condition (and it is henceforth finite-valued and continuous)

$\Delta_2$: $\exists c > 0, x_0 > 0$ such that $\forall x \geq x_0, \Psi(2x) \leq c\Psi(x) < +\infty$,

the two spaces $M^\Psi, L^\Psi$ coincide and $L^\Psi$ can be simply written as $\{f \in L^0 \mid E[\Psi(f)] < +\infty\} = \overline{L^\infty}^\Psi$ (where the closure is taken in the Luxemburg norm). This is the case of the $L^p$ spaces when $1 \leq p < +\infty$.

In [23], the authors also prove that when $\Psi$ is continuous on $\mathbb{R}$, then $M^\Psi = \overline{L^\infty}^\Psi$. So, when $\Psi$ is continuous, but grows too quickly, it may happen that $M^\Psi = \overline{L^\infty}^\Psi \subsetneq L^\Psi$. As a consequence, simple functions are not necessarily dense in $L^\Psi$ (see [23], Proposition III.4.3). This is quite a difference with classic $L^p$ spaces ($1 \leq p < +\infty$).

2.1.1. $L^{\widehat{u}}$ and $L^{\widehat{\Phi}}$. The even function $\widehat{u}\colon \mathbb{R} \to \mathbb{R} \cup \{+\infty\}$ defined by
$$\widehat{u}(x) = -u(-|x|) + u(0)$$
is a Young function and the induced Orlicz space is $L^{\widehat{u}} = \{f \in L^0 \mid \exists \alpha > 0,$ s.t. $E[\widehat{u}(\alpha f)] < +\infty\}$ with its Luxemburg norm $N_{\widehat{u}}(f)$.

REMARK 6. Note that $f \in L^{\widehat{u}}$ if and only if there exists some $\alpha > 0$ such that $E[u(-\alpha|f|)] > -\infty$. And $f \in M^{\widehat{u}}$ if and only if for all $\alpha > 0, E[u(-\alpha|f|)] > -\infty$.

As usual, the convex conjugate function $\widehat{\Phi}$ of $\widehat{u}$ is defined as
$$\widehat{\Phi}(y) \triangleq \sup_{x \in \mathbb{R}} \{xy - \widehat{u}(x)\}$$
and it is also a Young function. It admits a representation in terms of the convex conjugate $\Phi$ of the utility function $u$ as follows:
$$\widehat{\Phi}(y) = \begin{cases} 0, & \text{if } |y| \leq \beta, \\ \Phi(|y|) - \Phi(\beta), & \text{if } |y| > \beta, \end{cases}$$
where $\beta \geq 0$ is the right derivative of $\widehat{u}$ at 0, namely $\beta = D^+\widehat{u}(0) = D^-u(0)$, and $\Phi(\beta) = u(0)$. If $u$ is differentiable, note that $\beta = u'(0)$ and that $B$ is the unique solution of the equation $\Phi'(y) = 0$.

Let us consider the Orlicz space $L^{\widehat{\Phi}}$. It is convenient (see Section 2.2, item 2) to endow $L^{\widehat{\Phi}}$ with the Orlicz (or dual) norm
$$\|f\|_{\widehat{\Phi}} = \sup\{E[|fg|] \mid g \in L^{\widehat{u}} : E[\widehat{u}(g)] \leq 1\},$$
which is equivalent to the Luxemburg norm. As with all Orlicz spaces, $L^\infty \hookrightarrow L^{\widehat{\Phi}} \hookrightarrow L^1$.

REMARK 7. Obviously, $\Phi$ and $\widehat{\Phi}$ have the same behavior for large values, but $\widehat{\Phi}$ carries no information about the behavior of $\Phi$ near zero. For the comparison between $u$ and $\widehat{u}$, notice that $\widehat{u}(x)$ carries no information on the behavior of $u$ for $x > 0$, while for $x < 0$, we have simply $\widehat{u}(x) = -u(x) + u(0)$. This is a key point. In fact, when formulating the utility maximization problem in the Orlicz space $L^{\widehat{u}}$, we will only use this setting to control the losses of the terminal gains, that is, only the negative part of $(H \cdot S)_T$, or of the solution $f_x$, will belong to $L^{\widehat{u}}$.

*The case $a$ finite.* When the interior of the domain of $u$ is $(a, \infty)$ with $a < 0$ finite, evidently $\widehat{u}(x) = +\infty$ if $|x| > -a$. Since $u(-|x|) - u(0) \leq 0$ for all $x$, we have
$$xy - \widehat{u}(x) = xy + u(-|x|) - u(0) \leq -ay$$

for all $y > 0$ and $|x| \leq -a$ and therefore $\widehat{\Phi}(y) \leq -ay$. From these observations, $L^{\widehat{u}} = L^\infty, L^{\widehat{\Phi}} = L^1$ as sets and, trivially,

$$M^{\widehat{u}} = \{0\}.$$

Moreover, the identity map gives an isomorphism of Banach lattices, as the following lemma shows.

LEMMA 8. *The Luxemburg norm and the uniform norm on $L^{\widehat{u}} = L^\infty$ are equivalent. Consequently, also the Orlicz norm on $L^{\widehat{\Phi}} = L^1$ is equivalent to the $L^1$-norm.*

PROOF. Let $f \in L^{\widehat{u}} = L^\infty$. Then, if $E[\widehat{u}(\frac{f}{c})] \leq 1$, necessarily $\frac{|f|}{c} \leq -a$, so that

(9) $$\|f\|_\infty \leq -a N_{\widehat{u}}(f).$$

For the converse inequality, define $k$ to be the unique positive element of $(\widehat{u})^{-1}(\min(\widehat{u}(-a), 1))$. Evidently,

$$E\left[\widehat{u}\left(k\frac{f}{\|f\|_\infty}\right)\right] \leq 1,$$

whence

$$k N_{\widehat{u}}(f) \leq \|f\|_\infty. \qquad \square$$

When $a$ is finite, we recover the classical primal domain in the utility maximization problem (14): it is simply $(K^W - L^0_+) \cap L^\infty$ for the entire class of the utility functions with a half-line as proper domain.

Given Lemma 8, the explicit computations for $\widehat{\Phi}$ and $\widehat{u}$ are useless, but we give two examples for the sake of completeness.

1. Let $u$ be the logarithmic utility function and $a = -2$, that is,

$$u(x) = \begin{cases} \ln(2+x), & \text{if } x > -2, \\ -\infty, & \text{if } x \leq -2. \end{cases}$$

Then,

$$\widehat{u}(x) = \begin{cases} -\ln(2-|x|) + \ln 2, & \text{if } |x| < 2, \\ +\infty, & \text{if } |x| \geq 2, \end{cases}$$

while $\Phi(y) = -\ln y + 2y - 1$, $\beta = u'(0) = 1/2$, $\Phi(\beta) = u(0) = \ln 2$, so that

$$\widehat{\Phi}(y) = (\Phi(|y|) - \Phi(\beta)) I_{\{|y| > \beta\}} = (-\ln|y| + 2|y| - 1 - \ln 2) I_{\{|y| > 1/2\}}.$$

2. Suppose that $u(x) = \sqrt{4+x}$ if $x \geq -4$ and $u(x) = -\infty$ if $x < -4$. Then,

$$\widehat{u}(x) = \begin{cases} 2 - \sqrt{4 - |x|}, & \text{if } |x| \leq 4, \\ +\infty, & \text{if } |x| > 4, \end{cases}$$

and $\Phi(y) = \frac{1}{4y} + 4y$, so that

$$\widehat{\Phi}(y) = \left(\frac{1}{4|y|} + 4|y| - 2\right)I_{\{|y|>1/4\}}.$$

*The case $a = -\infty$.* Here, $\widehat{u}$ is continuous and, consequently, the subspace $M^{\widehat{u}} = \overline{L^{\infty}}^{\widehat{u}}$ is also a Banach space with the inherited $\widehat{u}$-norm. We give two examples, one with the exponential utility and the other with a linear utility.

1. When $u(x) = -e^{-x}$, $\widehat{u}(x) = e^{|x|} - 1$, while $\Phi(y) = y \ln y - y$ and $\widehat{\Phi}(y) = (|y| \ln |y| - |y| + 1)I_{\{|y| \geq 1\}}$. Therefore,

$$L^{\widehat{u}} = \{f \in L^0 \mid \exists \alpha > 0 \text{ s.t. } E[e^{\alpha|f|}] < +\infty\},$$

$$M^{\widehat{u}} = \{f \in L^0 \mid \forall \alpha > 0 \ E[e^{\alpha|f|}] < +\infty\}$$

and

$$L^{\widehat{\Phi}} = \{g \in L^0 \mid E[(|g| \ln |g|)I_{\{|g|>1\}}] < +\infty\}.$$

Due to convexity, we could remove the linear term from $\widehat{\Phi}$ in the above characterizations. Also, note that $M^{\widehat{u}}$ consists of those random variables that have *all* of the (absolute) exponential moments finite, while elements in $L^{\widehat{u}}$ are only required to have *some* finite exponential moments.

*In situations like the present one, the introduction of the Orlicz spaces shows its full potential.*

2. Let $u(x) = x$. Then, $\widehat{u}(x) = |x|$, $\Phi(y) = +\infty$ for all $y \geq 0$ and $\widehat{\Phi}(y) = (+\infty)I_{\{|y|>1\}} = \delta_{\{|y| \leq 1\}}$. So, $L^{\widehat{u}} = L^1 = M^{\widehat{u}}$ and $L^{\widehat{\Phi}} = L^{\infty}$.

In general, this is what the Orlicz spaces $L^{\widehat{u}}, L^{\widehat{\Phi}}$ reduce to whenever $u$ is asymptotically linear for $x \to -\infty$.

2.2. *On the norm dual of $L^{\widehat{u}}$ and $M^{\widehat{u}}$.* From the general theory of Banach lattices (see, e.g., [1]), we know that $(L^{\widehat{u}})^*$, the norm dual of $L^{\widehat{u}}$, admits the following decomposition:

$$(L^{\widehat{u}})^* = \mathcal{A} \oplus \mathcal{A}^d,$$

where $\mathcal{A}$ is the band of order-continuous linear functionals and $\mathcal{A}^d$ is the band of those singular ones which are *lattice orthogonal* to the functionals in $\mathcal{A}$. This means that every $z \in (L^{\widehat{u}})^*$ can be written in a unique way as $z = z_r + z_s$, with $z_r \in \mathcal{A}$ a regular functional, $z_s \in \mathcal{A}^d$ singular and $|z_r| \wedge |z_s| = 0$ [lattice orthogonality—we recall that the infimum $z^1 \wedge z^2$ of $z^1, z^2 \in (L^{\widehat{u}})^*$ can be characterized as the $z \in (L^{\widehat{u}})^*$ such that, for all positive $f \in L^{\widehat{u}}$, $z(f) = \inf_{0 \le g \le f}\{z^1(f-g) + z^2(g)\}$].

We can say more about the nature of the decomposition of $(L^{\widehat{u}})^*$. According to the specific nature of $u$, we distinguish between the two following cases.

1. If $a$ is finite, then $L^{\widehat{u}} = L^\infty$ and the above decomposition reduces to the Yosida–Hewitt one for elements of $ba(\Omega, \mathcal{F}, P)$,

$$ba = (L^\infty)^* = L^1 \oplus \mathcal{A}^d,$$

   where $\mathcal{A}^d$ consists of pure charges, that is, purely finitely additive measures.

2. If $a = -\infty$, then $\widehat{u}$ is continuous. For such Young functions, [3] and [21] showed that $\mathcal{A} = (M^{\widehat{u}})^*$ and that it can be identified with $L^{\widehat{\Phi}}$, endowed with the Orlicz (dual) norm. $\mathcal{A}^d$ is then the annihilator of $M^{\widehat{u}}$, denoted $(M^{\widehat{u}})^\perp$, whence

$$(L^{\widehat{u}})^* = (M^{\widehat{u}})^* \oplus (M^{\widehat{u}})^\perp = L^{\widehat{\Phi}} \oplus (M^{\widehat{u}})^\perp.$$

We remark that here, $M^{\widehat{u}} = \overline{L^\infty}^{\widehat{u}}$ and, consequently,

$$z \in \mathcal{A}^d \quad \text{if and only if} \quad \forall f \in L^\infty, z(f) = 0.$$

Therefore, we can identify the regular part $z_r \in \mathcal{A}$ of any $z \in (L^{\widehat{u}})^*_+$ with its density $\frac{dz_r}{dP} \in L^{\widehat{\Phi}}_+$ and we write its action on $f \in L^{\widehat{u}}$ as

$$z_r(f) = E_{z_r}[f].$$

REMARK 9. $L^{\widehat{u}} \subseteq L^1(Q)$ for all probabilities $Q$ such that $\frac{dQ}{dP} \in L^{\widehat{\Phi}}$. The space $L^{\widehat{\Phi}}$ can be identified with the regular elements in the dual of $L^{\widehat{u}}$, so this is a basic consequence of the general theory. However, we give here a simple and direct proof, which is based on the Fenchel inequality. Let us fix $f \in L^{\widehat{u}}$ and $\frac{dQ}{dP} \in L^{\widehat{\Phi}}$. Then, $E[\widehat{u}(\alpha f)] = E[\widehat{u}(\alpha|f|)] < \infty$ and $E[\widehat{\Phi}(\beta \frac{dQ}{dP})] < \infty$ for some positive $\alpha$ and $\beta$. From the Fenchel inequality $\alpha|f|\beta\frac{dQ}{dP} \le \widehat{u}(\alpha|f|) + \widehat{\Phi}(\beta\frac{dQ}{dP})$, we derive $f \in L^1(Q)$.

LEMMA 10. *The singular elements $z \in \mathcal{A}^d$ have the following characterization:*

(10) $\quad z \in \mathcal{A}^d \Leftrightarrow \forall f \in L^{\widehat{u}} \; \exists \text{ measurable } A_n \downarrow \varnothing \text{ such that } z(|f|I_{A_n}) = z(|f|).$

Proof. We may (and do) suppose that $z \geq 0$ (otherwise, we may work separately with $z^+, z^-$). The arrow ($\Leftarrow$) is immediate. In fact, the property $z(|f|I_{A_n}) = z(|f|)$ for some $A_n \downarrow \varnothing$ implies that $z \wedge \mu = 0$ for each regular $\mu$:

$$z \wedge \mu(f) = \inf_{0 \leq g \leq f} \{z(f-g) + \mu(g)\}$$

$$\leq \inf_n \{z(f - fI_{A_n}) + E_\mu[fI_{A_n}]\} = 0 \quad \text{for any } f \geq 0.$$

To prove ($\Rightarrow$), suppose that $f \geq 0$ and consider separately the two cases, $a$ finite and $a = -\infty$.

1. Case $a$ finite. Here, we can find a sequence $(A_n)_n$ which does not depend on the particular $f$. Since $z \in \mathcal{A}^d$, we have, in fact,

$$0 = z \wedge P(I_\Omega) = \inf_{0 \leq h \leq I_\Omega} \{z(I_\Omega - h) + E[h]\},$$

so there exists a sequence $0 \leq h^n \leq I_\Omega$ such that

$$0 \leq z(I_\Omega - h^n) + E[h^n] \leq \frac{1}{2^n}.$$

Call $\Omega_n = \{h^n > 0\}$. Then, the above inequalities, together with an application of the Borel–Cantelli lemma, imply that $\limsup_n \Omega_n = \varnothing$. So, by setting $A_n = \bigcup_{k \geq n} \Omega_k$, we have $A_n \downarrow \varnothing$ and

$$0 \leq z(I_\Omega - I_{A_n}) \leq z(I_\Omega - I_{\Omega_n}) \leq z(I_\Omega - h^n) \leq \frac{1}{2^n},$$

whence, necessarily, $z(I_\Omega - I_{A_n}) = 0$ for all $n$. This is equivalent to saying that $z$ is null on each $A_n^c$ and therefore for all $f \in L^{\widehat{u}} = L^\infty$, $z(f) = z(fI_{A_n})$.

2. Case $a = -\infty$. Take $A_n = \{f > n\}$ and consider the regular $\mu_n$ associated with $A_n$, namely $\mu_n(k) = E[kI_{A_n}]$. Since $z \in \mathcal{A}^d$, we have

$$0 = z \wedge \mu_n(f) = \inf_{0 \leq h \leq f} \{z(f - h) + E_{\mu_n}[h]\},$$

so there exists a sequence $0 \leq h^m \leq f$ such that

$$0 \leq z(f - h^m) + E_{\mu_n}[h^m] \leq \frac{1}{2^m}.$$

Therefore, for all $m$,

$$0 \leq z(f - fI_{A_n}) \leq z(f - h^m I_{A_n})$$

$$= z(f - h^m) + z(h^m I_{A_n^c}) \leq \frac{1}{2^m} + nz(I_{A_n^c}) = \frac{1}{2^m},$$

where the last equality holds since $z$ is null on $L^\infty$. Taking the limit over $m$, we obtain $z(f - fI_{A_n}) = 0$. □

The next proposition will be important in our applications, since it shows that $\mathcal{A}^d$ is an abstract Lebesgue space in the sense of [17]. That is, the norm of the dual space is additive on positive functionals in $\mathcal{A}^d$. We present a much simpler proof than that in [23], IV.3.4, based on Lemma 10.

PROPOSITION 11. *If $z \in \mathcal{A}^d_+$, then*

(11) $$\|z\| = \sup\{z(f) \mid f \geq 0, f \in L^{\widehat{u}} \text{ s.t. } E[\widehat{u}(f)] < +\infty\}.$$

*As a consequence, if $z_i \in \mathcal{A}^d_+$, then*

(12) $$\|z_1 + z_2\| = \|z_1\| + \|z_2\|.$$

PROOF. Call $l$ the supremum in (11). Since $z \geq 0$, we have $\|z\| = \sup\{z(f) \mid f \in B, f \geq 0\}$, where $B$ is the (open) unit ball of $L^{\widehat{u}}$. From the very definition of $B$, we see that $B \subseteq \{f \in L^{\widehat{u}} \text{ s.t. } E[\widehat{u}(f)] < +\infty\}$, so that $\|z\| \leq l$.

To show the opposite inequality, we use the characterization of $z$ provided in (10). Fix $f \in L^{\widehat{u}}_+$ so that $E[\widehat{u}(f)] < +\infty$. There then exists a sequence of measurable sets $A_n \downarrow \varnothing$ such that $z(fI_{A_n}) = z(f)$. But $fI_{A_n} \in B$ if $n$ is large enough. In fact, $\widehat{u}(fI_{A_n}) \downarrow 0$ and it is dominated by $\widehat{u}(f)$, so $E[\widehat{u}(fI_{A_n})]$ is definitely smaller than 1. We derive

$$l = \sup\{z(f) \mid f \geq 0, f \in L^{\widehat{u}} \text{ s.t. } E[\widehat{u}(f)] < +\infty\}$$
$$\leq \sup\{z(g) \mid g \geq 0, g \in B\} = \|z\|.$$

Additivity of the norm now follows easily, as in [23], Theorem 4.3.5. We sketch the proof. The only thing to show is that $\|z_1 + z_2\| \geq \|z_1\| + \|z_2\|$. This inequality can be obtained by taking positive functions $f_i \in L^{\widehat{u}}, i = 1, 2$, such that $E[\widehat{u}(f_i)] < +\infty$ and $z_i(f_i)$ is close to $\|z_i\|$, and observing that $f_1 \vee f_2 \in L^{\widehat{u}}$ and $E[\widehat{u}(f_1 \vee f_2)] < +\infty$. $\square$

**3. The utility maximization problem.** The conditions of compatibility and weak compatibility can now be expressed in the terminology of Orlicz space theory. In fact, a random variable $W \in L^0_+$ is:

- compatible iff $W \in M^{\widehat{u}}$;
- weakly compatible iff $W \in L^{\widehat{u}}$.

DEFINITION 12. When $W \in L^{\widehat{u}}_+$, the set of terminal values from admissible stochastic integrals is

$$K^W = \{(H \cdot S)_T \mid H \in \mathcal{H}^W\},$$

where $\mathcal{H}^W$ is defined in (4), and we set

(13) $$U^W(x) = \sup_{k \in K^W} E[u(x + k)].$$

REMARK 13.  Note that we do not require that $W$ is suitable because there is no need for this in constructing the duality. But, of course, this property is highly desirable. When $W \in L_+^{\widehat{u}} \cap \mathbb{S}$, the domain of maximization $K^W$ is nontrivial. What is more, under this stronger assumption, we will provide interesting characterizations of the dual variables (Proposition 19).

REMARK 14.  When $a$ is finite, $L^{\widehat{u}} = L^\infty$, so we could take $W = 1$ if we wished to recover the classical class $\mathcal{H}^1$ as the set of strategies with wealth bounded from below. This observation will be used in Section 6.1 for the comparison with [12].

In most of the preceding works on this subject, the basic idea for addressing the utility maximization problem (13) is to replace the domain $K^W$ with the set $(K^W - L_+^0) \cap L$, where $L$ is an appropriate topological vector space, for example, $L = L^\infty$, and to develop a dual approach based on the system $(L, L^*)$, where $L^*$ is the norm dual space of $L$.

LEMMA 15.  *Let $L$ be equal to either $L^\infty$, $M^{\widehat{u}}$ or $L^{\widehat{u}}$ and let $W \neq 0, W \in L_+$. If $g \in L_+$ and $k \in K^W$, then $k \wedge g \in (K^W - L_+^0) \cap L$. Moreover,*

$$(14) \qquad \sup_{k \in K^W} E[u(x+k)] = \sup_{f \in (K^W - L_+^0) \cap L} E[u(x+f)].$$

PROOF.  First, note that the hypothesis on $W$ excludes the possibility that $L = \{0\}$. So, either $L = L^\infty$ or, in case the utility is finite on the entire real line, $L = L^{\widehat{u}}$ or $L = M^{\widehat{u}}$. In all of these situations, $L \supseteq L^\infty$.

If $k \in K^W$, then there exists $c > 0$ such that $k^- \leq cW$, so that $k^- \in L_+$. Since $L$ is a Banach lattice, $k^+ \wedge g \in L_+$. Then, $k \wedge g = k^+ \wedge g - k^- \in L$ and it is also in $(K^W - L_+^0)$ since $k \wedge g = k - (k - k \wedge g) \in (K^W - L_+^0)$. To show (14), note that

$$\sup_{k \in K^W} E[u(x+k)] = \sup_{k \in K^W - L_+^0} E[u(x+k)]$$
$$\geq \sup_{f \in (K^W - L_+^0) \cap L^{\widehat{u}}} E[u(x+f)]$$
$$\geq \sup_{f \in (K^W - L_+^0) \cap L} E[u(x+f)].$$

To prove the other inequality, let $k \in K^W$ satisfy $E[u(x+k)] > -\infty$. Then, $E[u(-(x+k)^-)] > -\infty$, so for sufficiently large $n$, $E[u(-(x+k \wedge n)^-)] > -\infty$ and $E[u(x+k \wedge n)] > -\infty$. Hence, by monotone convergence, $E[u(x+k \wedge n)] \uparrow E[u(x+k)]$. The conclusion follows from $k \wedge n \in (K^W - L_+^0) \cap L$.  $\square$

The above lemma shows that the selection of the larger space $L = L^{\widehat{u}}$ is always consistent with the optimization over the set $K^W$ because the optimal values in (14) coincide. But it also ensures that whenever the behavior of the process $S$ is not too wild, so that not only $L^{\widehat{u}} \cap \mathbb{S}$, but also $M^{\widehat{u}} \cap \mathbb{S}$ (or even $L^\infty \cap \mathbb{S}$) is not empty, one could just as well use the smaller space $(K^W - L^0_+) \cap M^{\widehat{u}}$ [or $(K^W - L^0_+) \cap L^\infty$]. Set

$$C^W \triangleq (K^W - L^0_+) \cap L^{\widehat{u}}$$

and define $I_u : L^{\widehat{u}} \to [-\infty, \infty)$ by $I_u(f) = E[u(f)]$ and let $\mathcal{D}$ be the proper domain of $I_u$, that is,

$$\mathcal{D} \triangleq \{f \in L^{\widehat{u}} \mid E[u(f)] > -\infty\}.$$

PROPOSITION 16. *The concave functional $I_u$ on $L^{\widehat{u}}$ is proper and it is norm-continuous on the interior of its proper domain, which is not empty. Moreover, there exists a norm continuity point of $I_u$ that belongs to $C^W$.*

PROOF. Thanks to [15], Proposition I.2.5, the thesis is equivalent to showing that there is a nonempty open set $O$ on which $I_u$ is not everywhere equal to $+\infty$ and bounded below by a constant $c \in \mathbb{R}$. We show slightly more, that is, on the open unit ball $B$ of $L^{\widehat{u}}$, the functional $I_u$ is (i) everywhere less than $+\infty$ and (ii) uniformly bounded below.

(i) If $b \in B$, then $E[|b|] < +\infty$, so, by Jensen's inequality, $I_u(b) \le u(E[b]) < +\infty$.

(ii) For all $b \in B$, $E[\widehat{u}(b)] \le 1$ and $E[\widehat{u}(b^-)] \le 1$. Hence,

$$-I_u(-b^-) = E[-u(-b^-)] = E[\widehat{u}(b^-)] - u(0) \le 1 - u(0)$$

and so $I_u(b) \ge I_u(-b^-) \ge u(0) - 1$. Note for future use that (i) and (ii) clearly imply that $I_u$ is finite on the ball $B$.

The second statement of the proposition follows from $C^W \supseteq -B_+$. □

The next lemma is a very nice consequence of the choice of the right Orlicz space $L^{\widehat{u}}$.

LEMMA 17. *Let $z \in \mathcal{A}^d_+$. Then,*

(15) $$\|z\| = \sup_{f \in \mathcal{D}} z(-f).$$

*In the case $a$ finite, $z$ is a nonnegative pure charge and*

(16) $$\|z\| = -az(\Omega).$$

PROOF. Since $z \geq 0$, $\sup_{f \in \mathcal{D}} z(-f) = \sup_{f \leq 0, f \in \mathcal{D}} z(-f)$. But $f \leq 0$, $f \in \mathcal{D}$ if and only if $g = -f$ is a nonnegative random variable satisfying $E[\widehat{u}(g)] < +\infty$. The thesis then follows from (11). The case $a$ finite is then obvious: $\|z\| = \sup_{\{f \in L_+^{\widehat{u}}, f < -a\}} z(f) = -az(\Omega)$. □

3.1. *Dual variables.* We now describe the dual variables. In what follows, $W \in L_+^{\widehat{u}}$, and we always refer to the dual system $(L^{\widehat{u}}, (L^{\widehat{u}})^*)$. Consider the polar cone

$$(C^W)^0 = \{z \in (L^{\widehat{u}})^* \mid z(f) \leq 0 \ \forall f \in C^W\}$$

and define

(17)  $$\mathcal{M}^W \triangleq \{Q \in (C^W)^0 \mid Q(I_\Omega) = 1\}.$$

Note that $(C^W)^0 \subseteq (L^{\widehat{u}})_+^*$ since $(-L_+^{\widehat{u}}) \subseteq C^W$. Therefore, the functionals of interest are positive. We also remark that in the case $a = -\infty$, the condition in (17) amounts to saying that $E_{Q_r}[I_\Omega] = 1$ since $Q_s$ vanishes over $L^\infty$. So, if $Q \in \mathcal{M}^W$ is regular, then $\frac{dQ}{dP}$ is a probability density. In Proposition 19 below, we completely characterize the absolutely continuous probability measures arising in this way.

REMARK 18. Notice that the regular elements in $\mathcal{M}^W$ can be described as

(18)  $$\mathcal{M}^W \cap L^{\widehat{\Phi}} = \left\{Q \ll P \ \Big| \ \frac{dQ}{dP} \in L^{\widehat{\Phi}} \text{ and } E_Q[f] \leq 0 \ \forall f \in C^W\right\}$$

(19)  $$= \left\{Q \ll P \ \Big| \ \frac{dQ}{dP} \in L^{\widehat{\Phi}} \text{ and } E_Q[(H \cdot S)_T] \leq 0 \ \forall H \in \mathcal{H}^W\right\}.$$

Indeed, the set in (19) is clearly contained in (18). To check the opposite inclusion, let $Q \in \mathcal{M}^W \cap L^{\widehat{\Phi}}$, $H \in \mathcal{H}^W$ and note that $E_Q[(H \cdot S)_T]$ is well defined since $(H \cdot S)_T \geq -cW$ and $W \in L^1(Q)$ for all $Q \in L^{\widehat{\Phi}}$ (Remark 9). Furthermore, from Lemma 15, $(H \cdot S)_T \wedge n \in C^W$ for each $n$ and, by monotone convergence, $E_Q[(H \cdot S)_T] = \lim_n E_Q[(H \cdot S)_T \wedge n] \leq 0$.

We denote by $\mathbb{M}_\sigma$ the set of $P$-absolutely continuous $\sigma$-martingale probabilities for $S$ (see [13, 14] for more information about this concept). Recall that when $S$ is bounded (resp. locally bounded), we have

$$\mathbb{M}_\sigma = \{Q \ll P \mid S \text{ is a martingale (resp. local martingale) w.r.t. } Q\},$$

that is, $\mathbb{M}_\sigma$ is the set of $P$-a.c. *martingale* (resp. *local martingale*) probabilities.

Set
$$\mathbb{M}_{\sup}^W = \{Q \ll P \mid H \cdot S \text{ is a } Q\text{-supermartingale } \forall H \in \mathcal{H}^W\},$$
$$\mathcal{H}^{\widehat{u}} = \bigcup_{W \in L^{\widehat{u}}, W \geq 1} \mathcal{H}^W.$$

PROPOSITION 19.  *The regular elements in $\mathcal{M}^W$ have the following, interesting, probabilistic properties:*

(a) *if $W \in L_+^{\widehat{u}}$, then*
$$\mathbb{M}_\sigma \cap L^{\widehat{\Phi}} \subseteq \mathbb{M}_{\sup}^W \cap L^{\widehat{\Phi}};$$

(b)
$$\mathbb{M}_\sigma \cap L^{\widehat{\Phi}} \subseteq \{Q \ll P \mid Q \in L^{\widehat{\Phi}} \text{ and } H \cdot S \text{ is a } Q\text{-supermartingale } \forall H \in \mathcal{H}^{\widehat{u}}\};$$

(c) *if $W \in L^{\widehat{u}}$, $W \geq 1$, then*
$$\mathbb{M}_{\sup}^W \cap L^{\widehat{\Phi}} = \mathcal{M}^W \cap L^{\widehat{\Phi}};$$

(d) *if $W \in L^{\widehat{u}} \cap \mathbb{S}$, then*

(20) $$\mathbb{M}_\sigma \cap L^{\widehat{\Phi}} = \mathbb{M}_{\sup}^W \cap L^{\widehat{\Phi}} = \mathcal{M}^W \cap L^{\widehat{\Phi}};$$

(e) *if $L^{\widehat{u}} \cap \mathbb{S} \neq \varnothing$, then*
$$\mathbb{M}_\sigma \cap L^{\widehat{\Phi}} = \{Q \ll P \mid Q \in L^{\widehat{\Phi}} \text{ and } H \cdot S \text{ is a } Q\text{-supermartingale } \forall H \in \mathcal{H}^{\widehat{u}}\}.$$

PROOF. (a) Let $Q \in \mathbb{M}_\sigma \cap L^{\widehat{\Phi}}$ and $W \in L_+^{\widehat{u}}$. Since $Q \in L^{\widehat{\Phi}}$, $W \in L^1(Q)$. If $H \in \mathcal{H}^W$, then there exists a $c \geq 0$ such that $(H \cdot S)_t \geq -cW$. From $Q \in \mathbb{M}_\sigma$, we can find a positive predictable scalar process $\psi$ such that $\psi^{-1} \cdot S^i$ is a $Q$ uniformly integrable martingale for $i = 1, \ldots, d$.

If we set $X$ to be the semimartingale with $i$th component $X^i = \psi^{-1} \cdot S^i$, we can write $H \cdot S = (\psi H) \cdot X$. That is, $H \cdot S$ is a stochastic integral with respect to the $Q$-martingale $X$ and its negative part is controlled by the $Q$-integrable variable $cW$. Thanks to a lemma of Ansel and Stricker [4], $H \cdot S$ is then a $Q$-local martingale and a supermartingale. Hence, $Q \in \mathbb{M}_{\sup}^W$.

(b) This follows from (a) and the fact that

(21) $$\bigcap_{W \in L^{\widehat{u}}, W \geq 1} \mathbb{M}_{\sup}^W = \{Q \ll P \mid H \cdot S \text{ is a } Q\text{-supermartingale } \forall H \in \mathcal{H}^{\widehat{u}}\}.$$

(c) Obviously, $\mathbb{M}_{\sup}^W \cap L^{\widehat{\Phi}} \subseteq \mathcal{M}^W \cap L^{\widehat{\Phi}}$ and so it remains to show that $\mathcal{M}^W \cap L^{\widehat{\Phi}} \subseteq \mathbb{M}_{\sup}^W \cap L^{\widehat{\Phi}}$. Define the stopping times (increasing to $T$) $T_n = \inf\{t \leq T | (H \cdot S)_t > n\}$ and fix $s < t \leq T$ and $A \in F_s$. Let $W \in L^{\widehat{u}}$ and $W \geq 1$. If $(H \cdot S)_t \geq -cW$, then also $I_A I_{]s,t \wedge T_n]} H \in \mathcal{H}^W$ since

$$((I_A H I_{]s,t \wedge T_n]}) \cdot S)_u \geq -cW - n \geq -(c+n)W.$$

When $Q \in \mathcal{M}^W \cap L^{\widehat{\Phi}}$, we have $E_Q[((I_A H I_{]s,t \wedge T_n]}) \cdot S)_T] \leq 0$ so that

$$E_Q[I_A (H \cdot S)_{t \wedge T_n} I_{\{T_n > s\}}] \leq E_Q[I_A (H \cdot S)_s I_{\{T_n > s\}}].$$

Now, observe that $|I_A (H \cdot S)_s I_{\{T_n > s\}}| \leq |(H \cdot S)_s|$ and $(H \cdot S)_s \in L^1(Q)$, since $I_{[0,s]} H \in \mathcal{H}^W$. In addition, $I_A (H \cdot S)_{t \wedge T_n} I_{\{T_n > s\}} \geq -cW$, hence an application of the Lebesgue dominated convergence theorem on the RHS and of Fatou's lemma on the LHS leads us to the desired inequality $E_Q[I_A (H \cdot S)_t] \leq E_Q[I_A (H \cdot S)_s]$.

(d) We only need to show that $\mathcal{M}^W \cap L^{\widehat{\Phi}} \subseteq \mathbb{M}_\sigma \cap L^{\widehat{\Phi}}$. Suppose that $Q \in \mathcal{M}^W \cap L^{\widehat{\Phi}} \subseteq L_+^{\widehat{\Phi}}$. It is easily checked that the random variables $\pm (H^i I_A I_{]s,t]} \cdot S^i)_T$ satisfy

$$-2W \leq \pm(H^i I_A I_{]s,t]} \cdot S^i)_T \leq 2W$$

for all $s < t$, $A \in \mathcal{F}_s$, where the integrands $H^i$ are those in the definition of suitable $W$, relation (5). If we let $k_T = (H^i I_A I_{]s,t]} \cdot S^i)_T$, then $\pm k_T \in K^W \cap L^{\widehat{u}}$, so, by definition of $(C^W)^0$, we deduce $E_Q[k_T] = 0$. Hence, for all $i = 1, \ldots, d$, $H^i \cdot S^i$ is a $Q$-martingale. This implies that $S^i$ is a $\sigma$-martingale with respect to $Q$, thanks to the characterization provided by [14].

(e) Note that

$$\mathbb{M}_\sigma \cap L^{\widehat{\Phi}} \subseteq \bigcap_{W \in L^{\widehat{u}}, W \geq 1} \mathbb{M}_{\sup}^W \cap L^{\widehat{\Phi}} \subseteq \bigcap_{W \in L^{\widehat{u}} \cap \mathbb{S}} \mathbb{M}_{\sup}^W \cap L^{\widehat{\Phi}} = \mathbb{M}_\sigma \cap L^{\widehat{\Phi}},$$

where the first inclusion follows from (a) and the equality from (d). The thesis then follows from (21). □

COROLLARY 20. *If $L^{\widehat{u}} \cap \mathbb{S} \neq \varnothing$, $\mathcal{M}^W \cap L^{\widehat{\Phi}}$ does not depend on which $W$ is selected in $L^{\widehat{u}} \cap \mathbb{S}$ and it coincides with $\mathbb{M}_\sigma \cap L^{\widehat{\Phi}}$.*

3.2. *Minimax theorem.*

THEOREM 21. *Let $u : \mathbb{R} \to \mathbb{R} \cup \{-\infty\}$ be increasing and concave on the interior $(a, \infty)$, $a \in [-\infty, 0)$, of its effective domain and with the property $\lim_{x \to -\infty} u(x) = -\infty$.*

*If there exists $W \in L_+^{\widehat{u}}$ satisfying $\sup_{k \in K^W} E[u(x+k)] < u(+\infty)$ for some $x > a$, then $\mathcal{M}^W$ is not empty and*

$$(22) \qquad U^W(x) \triangleq \sup_{k \in K^W} E[u(x+k)] = \sup_{f \in C^W} E[u(x+f)]$$

$$(23) \qquad \qquad = \min_{\lambda > 0, Q \in \mathcal{M}^W} \left\{ \lambda(x + \|Q_s\|) + E\left[\Phi\left(\lambda \frac{dQ_r}{dP}\right)\right] \right\},$$

*where $Q = Q_r + Q_s$ is the decomposition of $Q$ into regular and singular part.*

PROOF. We first prove the result in the case $x = 0$. The concave conjugate functional of $I_u$ is $J_u(z) = -E[\Phi(\frac{dz_r}{dP})] - \sup_{f \in \mathcal{D}} z_s(-f)$ by [18], Theorem 2.6 (see Theorem 46 in the Appendix). From (15), $J_u(z) = -E[\Phi(\frac{dz_r}{dP})] - \|z_s\|$. By Proposition 16, the Fenchel duality theorem can be applied (see, e.g., Brezis [11]) to obtain

$$\sup_{f \in C^W} E[u(f)] = \min_{z \in (C^W)^0} -J_u(z)$$

$$(24) \qquad \qquad = \min_{z \in (C^W)^0} \left\{ E\left[\Phi\left(\frac{dz_r}{dP}\right)\right] + \|z_s\| \right\}$$

$$(25) \qquad \qquad = \min_{\lambda > 0, Q \in \mathcal{M}^W} \left\{ E\left[\Phi\left(\lambda \frac{dQ_r}{dP}\right)\right] + \lambda \|Q_s\| \right\},$$

where the last equality follows from a reparametrization via $\mathcal{M}^W$. This last passage is ensured by the assumption $\sup_{k \in K^W} E[u(k)] < u(+\infty)$, together with $\Phi(0) = u(\infty)$. From these conditions, in fact, we derive that any solution $z^*$ of the dual problem in (24) has nonzero regular component, so that $\mathcal{M}^W \neq \varnothing$, even when $a = -\infty$.

The case with arbitrary initial endowment $x$ follows from a few considerations. If we let $u_x(\cdot) = u(x + \cdot)$, then $u_x$ is finite on $(a_x, \infty)$, $a_x = a - x < 0$, $I_{u_x}(f) = E[u(x+f)]$ and the proper domain of $I_{u_x}$ is $\mathcal{D}_x = \mathcal{D} - x$. Then, we go on as in case $x = 0$, taking into account that the concave conjugate of $I_{u_x}$ is $J_{u_x}(z) = -xz_r(\Omega) - E[\Phi(\frac{dz_r}{dP})] - \sup_{f \in \mathcal{D}_x} z_s(-f)$. And, since $\sup_{f \in \mathcal{D}_x} z_s(-f) = \sup_{g \in \mathcal{D}} z_s(-g) + xz_s(\Omega) = \|z_s\| + xz_s(\Omega)$, we conclude that $J_{u_x}(z) = -xz(\Omega) - E[\Phi(\frac{dz_r}{dP})] - \|z_s\|$. □

REMARK 22. (i) Suppose that $u$ is strictly concave so that $E[\Phi(\cdot)]$ is strictly convex. The optimal functional $Q^*$ is then unique only in the regular part $Q_r^*$. In fact, $\|\cdot\|$ is additive on the nonnegative singular functionals (Proposition 11). Therefore, the dual objective function to be minimized in (24) is not strictly convex and the nonuniqueness of the solution can only arise from the singular part.

(ii) Even if the duality is formulated with respect to the dual system $(L^{\widehat{u}}, L^{\widehat{\Phi}} \oplus \mathcal{A}^d)$, we have to keep in mind that it is the function $\Phi$, the conjugate of $u$, that shows up in the dual problem (25), and not $\widehat{\Phi}$. This is also the reason why, in the next section, we have to consider the smaller set $\mathcal{N}_\Phi^W$ instead of $\mathcal{M}^W$.

Until now, we have shown that there is no duality gap between the primal problem (22) and the dual problem (23) and that the dual problem is attained.

**4. Dual and primal optima.** In this section, we analyze the properties of the solutions of the dual problem (23) and prove the existence of the solution to the primal problem over a set larger than $K^W$, which is defined in (33) below. As should be clear from Remark 7, we may not expect that the optimal solution $f_x$ belongs to $L^{\widehat{u}}$, but only that $f_x^- \in L^{\widehat{u}}$. For example, think of the case $a$ finite. Then, $L^{\widehat{u}} = L^\infty$ and it is well known that the primal solution may not be bounded.

*In this section and in the sequel, we will work under assumptions* (A1) *and* (A2).

The convex conjugate $\Phi$ is then a strictly convex differentiable function satisfying $\Phi(+\infty) = +\infty$, $\Phi(0^+) = u(+\infty)$, $\Phi' = -(u')^{-1}$, $\Phi'(0^+) = -\infty$ and $\Phi'(+\infty) = -a$.

Let

$$\mathcal{L}_\Phi \triangleq \{f \in L_+^0(P) \mid E[\Phi(f)] < \infty\} \subseteq L_+^{\widehat{\Phi}}$$

*and note that assumption* (A2) *is equivalent to requiring that $\mathcal{L}_\Phi$ is a convex cone.*

REMARK 23. While, by construction, $\widehat{\Phi}(0) = 0$, $\Phi(0) = \infty$ is possible. So, in general, we have only $\mathcal{L}_\Phi \subseteq L_+^{\widehat{\Phi}}$. Of course, $\mathcal{L}_\Phi = L_+^{\widehat{\Phi}}$ whenever $\Phi(0) < \infty$ [which is equivalent to $u(+\infty) < +\infty$].

REMARK 24 [*Consequences of* (A2)]. Assumption (A2) implies that $L^{\widehat{\Phi}} = M^{\widehat{\Phi}}$. Indeed, it implies that if $f \in L_+^0(P)$, then

$$E[\Phi(f)\mathbf{1}_{\{f \geq \beta\}}] < +\infty \Leftrightarrow E[\Phi(\lambda f)\mathbf{1}_{\{f \geq \beta\}}] < +\infty \quad \forall \lambda > 0,$$

and therefore

$$E[\widehat{\Phi}(f)] < +\infty \quad \text{iff} \quad E[\widehat{\Phi}(\lambda f)] < +\infty \quad \forall \lambda > 0.$$

Due to (A2), when $E[\Phi(\lambda \frac{dQ_r}{dP})] < \infty$, we have $Q_r \in \mathcal{L}_\Phi \subseteq L_+^{\widehat{\Phi}}$. Therefore, the min in (23) is reached on the convex set of functionals

$$\mathcal{N}_\Phi^W \triangleq \mathcal{M}^W \cap \{Q \in (L^{\widehat{u}})^* \mid Q_r \neq 0, Q_r \in \mathcal{L}_\Phi\}.$$

We will simply write $\mathcal{N}^W$ instead of $\mathcal{N}^W_\Phi$.

The first two propositions are extensions of some of the results in [5], where only probability measures $Q$ were allowed.

PROPOSITION 25. *Suppose that $\mathcal{N} \subseteq (L^{\widehat{u}})^*_+$ is a convex set and that, for each $Q = Q_r + Q_s \in \mathcal{N}$, $E_{Q_r}[I_\Omega] > 0$. Let $\mathcal{N}_r = \{Q_r \mid Q \in \mathcal{N}\}$ and suppose that $\mathcal{N}_r \subseteq \mathcal{L}_\Phi$. If $Q^\lambda \in \mathcal{N}$ is optimal for*

$$\text{(26)} \qquad \inf\left\{ E\left[\Phi\left(\lambda \frac{dQ_r}{dP}\right)\right] + \lambda \|Q_s\| \mid Q \in \mathcal{N} \right\},$$

*then, $\forall Q \in \mathcal{N}$*

$$\text{(27)} \qquad E_{Q_r}\left[-\Phi'\left(\lambda \frac{dQ_r^\lambda}{dP}\right)\right] - \|Q_s\| \leq E_{Q_r^\lambda}\left[-\Phi'\left(\lambda \frac{dQ_r^\lambda}{dP}\right)\right] - \|Q_s^\lambda\|.$$

PROOF. First consider $\lambda = 1$. Suppose that $Q^1 = Q_r^1 + Q_s^1$ is optimal for (26) and let $Q^0 = Q_r^0 + Q_s^0 \in \mathcal{N}$. Set $\eta^0 = \frac{dQ_r^0}{dP}$, $\eta^1 = \frac{dQ_r^1}{dP}$, $\eta^x = x\eta^1 + (1-x)\eta^0$, $Q_s^x = xQ_s^1 + (1-x)Q_s^0$, $x \in [0, 1]$. From the convexity of $\Phi$, we derive

$$\text{(28)} \qquad \eta^0 \Phi'(\eta^1) \leq \eta^1 \Phi'(\eta^1) + \Phi(\eta^0) - \Phi(\eta^1), \qquad P\text{-a.s.}$$

As in Lemma 2 of [5], we again exploit the convexity of $\Phi$, the hypothesis $\mathcal{N}_r \subseteq \mathcal{L}_\Phi$ and the cone property of $\mathcal{L}_\Phi$ guaranteed by (A2) to deduce from (28) the two integrability conditions

$$\eta^1 \Phi'(\eta^1) \in L^1(P) \quad \text{and} \quad (\eta^0 \Phi'(\eta^1))^+ \in L^1(P).$$

Set $F(x) = \Phi(\eta^x)$, $x \in [0, 1]$. By convexity of $F$, $(\frac{F(1)-F(x)}{1-x})$ is monotone. Since $E[F(1) - F(0)]$ is finite, we apply the monotone convergence theorem to obtain

$$\left(\frac{d}{dx} E[\Phi(\eta^x)]\Big|_{x=1}\right) \triangleq \lim_{x \uparrow 1} E\left[\frac{F(x) - F(1)}{x-1}\right]$$
$$= E[F'(1)] = E[\Phi'(\eta^1)(\eta^1 - \eta^0)].$$

From $(\eta^x\, dP + Q_s^x) \in \mathcal{N}$, the linearity of $x \to \|Q_s^x\|$ and optimality of $Q^1$, we see that the left derivative at $x = 1$ is negative:

$$\text{(29)} \qquad \frac{d}{dx}\{E[\Phi(\eta^x)] + \|Q_s^x\|\}\Big|_{x=1} \leq 0.$$

So

$$E[\Phi'(\eta^1)(\eta^1 - \eta^0)] + \|Q_s^1\| - \|Q_s^0\| \leq 0$$

and we have (27).

Again due to the cone property of $\mathcal{L}_\Phi$, the case with general $\lambda > 0$ follows from the case with $\lambda = 1$, applied to the functions $\Phi_\lambda(y) = \Phi(\lambda y)$ and $\lambda\|\cdot\|$. □

PROPOSITION 26. *If $Q \in \mathcal{L}_\Phi$ and $0 < Q(\Omega) \leq 1$, then, for all $c > aQ(\Omega)$, the optimal $\lambda(c; Q)$ solution of*

$$\min_{\lambda > 0} \left\{ \lambda c + E\left[\Phi\left(\lambda \frac{dQ}{dP}\right)\right] \right\} \tag{30}$$

*is the unique positive solution of the first-order condition*

$$c + E\left[\frac{dQ}{dP}\Phi'\left(\lambda \frac{dQ}{dP}\right)\right] = 0. \tag{31}$$

*The random variable $f^* \triangleq -\Phi'(\lambda(c;Q)\frac{dQ}{dP}) \in \{f \in L^1(Q) \mid E_Q[f] = c\}$ satisfies $u(f^*) \in L^1(P)$ and*

$$\min_{\lambda > 0} \left\{ \lambda c + E\left[\Phi\left(\lambda \frac{dQ}{dP}\right)\right] \right\} = \sup\{E[u(f)] \mid f \in L^1(Q) \text{ and } E_Q[f] \leq c\} \tag{32}$$
$$= E[u(f^*)] < u(\infty).$$

*In the case where $a$ is finite, if $c \leq aQ(\Omega)$, then the min in (30) is not obtained and the optimal value is $-\infty$.*

PROOF. It follows, as in the proof of [5], Proposition 7, replacing $x$ with $c$ and noting (see also [5], Lemma 2-c) that now, $F(\lambda) = E[\frac{dQ}{dP}\Phi'(\lambda \frac{dQ}{dP})]$ is a bijection between $(0, +\infty)$ and $(-\infty, -aQ(\Omega))$. □

Under the same hypothesis and notation of Theorem 21 and Remark 24, we know that $\mathcal{N}^W \neq \varnothing$, so we may define

$$K_\Phi^W(x) = \{f \in L^0 \mid f \in L^1(Q_r), E_{Q_r}[f] \leq x + \|Q_s\| \ \forall Q \in \mathcal{N}^W\} \tag{33}$$

and

$$U_\Phi^W(x) = \sup_{f \in K_\Phi^W(x)} E[u(f)].$$

LEMMA 27. *Let $W \in L_+^{\widehat{u}}$. If $k \in K^W$ and $E[u(x+k)] > -\infty$, then $x + k \in K_\Phi^W(x)$. Moreover,*

$$U^W(x) \leq U_\Phi^W(x) \leq \inf_{\lambda > 0, Q \in \mathcal{N}^W} \left\{ \lambda(x + \|Q_s\|) + E\left[\Phi\left(\lambda \frac{dQ_r}{dP}\right)\right] \right\} \tag{34}$$

*and, under the same hypothesis of Theorem 21, the inf is attained.*

PROOF. Since $k \in K^W$, by Lemma 15, we have $k \wedge n \in C^W$. Hence, $Q(x + k \wedge n) = Q(x) + Q(k \wedge n) \leq x$ for all $Q \in \mathcal{N}^W$. From this, $E_{Q_r}[x +$ $k \wedge n] \leq x - Q_s(x + k \wedge n)$. If $n > -x$, we have $(x + k \wedge n)^- = (x + k)^-$. By assumption, $-(x+k)^- \in \mathcal{D}$ so that, for sufficiently large $n$,

$$(35) \qquad E_{Q_r}[x + k \wedge n] \leq x - Q_s(x + k \wedge n) \leq x + Q_s((x + k \wedge n)^-)$$

$$(36) \qquad\qquad\qquad = x + Q_s((x+k)^-) \leq x + \|Q_s\|,$$

where the last inequality follows from (15). We then conclude that $x + k \in K_\Phi^W(x)$ and $U^W(x) \leq U_\Phi^W(x)$. From the Fenchel inequality, we have

$$u(f) \leq \lambda \frac{dQ_r}{dP} f + \Phi\left(\lambda \frac{dQ_r}{dP}\right) \qquad \forall f \in K_\Phi^W(x), \forall Q \in \mathcal{N}^W, \forall \lambda > 0.$$

By taking the expectations and optimizing on both sides, we obtain the second inequality in (34). $\square$

LEMMA 28. *Let $W \in L_+^{\widehat{u}}$. If $\mathcal{N}^W \neq \varnothing$, then $U^W(x) < u(\infty)$ for all $x > a$.*

PROOF. If $Q \in \mathcal{N}^W$, then, from Lemma 27, we deduce

$$U^W(x) \leq \sup_{f \in K_\Phi^W(x)} E[u(f)]$$

$$(37) \qquad\qquad \leq \sup\{E[u(f)] \mid f \in L^1(Q_r) \text{ and } E_{Q_r}[f] \leq x + \|Q_s\|\}.$$

Since $x > a$, we have $x + \|Q_s\| > aQ_r(\Omega)$ [in the case $a = -\infty$, this is certainly true because $Q_r(\Omega) = 1 \neq 0$; in the case $a$ finite, it follows from $\|Q_s\| = -aQ_s(\Omega)$, as shown in (16)]. Hence, an application of Proposition 26, (32), gives that the last term in (37) is less than $u(\infty)$. $\square$

We are ready to state and prove the main result of this section. As we have repeated throughout the paper, the primal solution does not generally belong to the set $K^W$. In [5], we showed that in the case $a = -\infty$, this can happen, even if $S$ is locally bounded.

THEOREM 29. *Assume that* (A1) *and* (A2) *hold. If there exists $W \in L_+^{\widehat{u}}$ that satisfies*

$$(38) \qquad \sup_{k \in K^W} E[u(x + k)] < u(\infty) \qquad \text{for some } x > a,$$

*then $\mathcal{N}^W$ is not empty and, for all $x > a$:*

1.

$$U^W(x) = U_\Phi^W(x) = E[u(f_x)]$$

$$(39) \qquad\qquad = \min_{\lambda > 0, Q \in \mathcal{N}^W}\left\{\lambda(x + \|Q_s\|) + E\left[\Phi\left(\lambda \frac{dQ_r}{dP}\right)\right]\right\};$$

2. *there exists a unique solution to*

$$U_\Phi^W(x) = \sup\{E[u(f)] \mid f \in K_\Phi^W(x)\}$$

*and it is given by*

(40) $$f_x \triangleq -\Phi'\left(\lambda^* \frac{dQ_r^*}{dP}\right) \in K_\Phi^W(x),$$

*where $\lambda^*$ (unique) and $Q^*$ (unique in the regular part) are solutions to the dual problem in (39);*

3.

(41) $$E_{Q_r^*}[f_x] = x + \|Q_s^*\|.$$

PROOF. From Theorem 21, Remark 24 and Lemma 28, $\mathcal{N}^W$ is not empty and $U^W(x) < u(\infty)$ for all $x > a$. Therefore, from Theorem 21, Lemma 27 and Remark 24,

$$U^W(x) = U_\Phi^W(x) = \min_{\lambda > 0, Q \in \mathcal{N}^W}\left\{\lambda(x + \|Q_s\|) + E\left[\Phi\left(\lambda \frac{dQ_r}{dP}\right)\right]\right\}$$
$$= \lambda^*(x + \|Q_s^*\|) + E\left[\Phi\left(\lambda^* \frac{dQ_r^*}{dP}\right)\right].$$

As in the proof of Lemma 28, from $x > a$, we deduce $x + \|Q_s^*\| > aQ_r^*(\Omega)$. Hence, from Proposition 26, we obtain that $\lambda^* = \lambda^*(x + \|Q_s^*\|; Q_r^*)$ is the unique solution to

$$x + \|Q_s^*\| + E\left[\frac{dQ_r^*}{dP}\Phi'\left(\lambda \frac{dQ_r^*}{dP}\right)\right] = 0$$

so that the r.v. $f_x \triangleq -\Phi'(\lambda^* \frac{dQ_r^*}{dP})$ satisfies $E_{Q_r^*}[f_x] = x + \|Q_s^*\|$ and $E[u(f_x)] = \lambda^*(x + \|Q_s^*\|) + E[\Phi(\lambda^* \frac{dQ_r^*}{dP})]$. In addition, since $Q^*$ is an optimal solution of

$$\min_{Q \in \mathcal{N}^W}\left\{E\left[\Phi\left(\lambda^* \frac{dQ_r}{dP}\right)\right] + \lambda^*\|Q_s\|\right\},$$

we can apply Proposition 25 to deduce

$$E_{Q_r}\left[-\Phi'\left(\lambda^* \frac{dQ_r^*}{dP}\right)\right] - \|Q_s\| \leq E_{Q_r^*}\left[-\Phi'\left(\lambda^* \frac{dQ_r^*}{dP}\right)\right] - \|Q_s^*\| = x,$$

which means that $f_x \in K_\Phi^W(x)$. □

4.1. *Further properties and comments.* Since the optimum satisfies $E[u(f_x)] > -\infty$, it follows that $f_x^- \in L^{\widehat{u}}$ and, from Lemma 17, that $\|Q_s^*\| = \sup_{f \in \mathcal{D}} Q_s^*(-f) \geq Q_s^*(f_x^-)$. We now provide a sufficient condition to obtain the equality $\|Q_s^*\| = Q_s^*(f_x^-)$. First, we need a simple lemma.

LEMMA 30. *Let $f \leq 0$, $f \in \mathcal{D}$. There then exists an $\epsilon > 0$ such that $(1+\epsilon)f \in \mathcal{D}$ iff $f \in int(\mathcal{D})$.*

PROOF. We prove only one implication, the other being trivial. If $(1+\epsilon)f \in \mathcal{D}$, then $f + \frac{\epsilon}{1+\epsilon}B \subset \mathcal{D}$, where $B$ is the unit open ball of $L^{\widehat{u}}$. In fact, if $g \in B$, then, by the proof of Proposition 16, $E[u(g)]$ is finite and $E[u(f + \frac{\epsilon}{1+\epsilon}g)] \geq \frac{E[u((1+\epsilon)f)] + \epsilon E[u(g)]}{1+\epsilon} > -\infty$. □

PROPOSITION 31. *Suppose that the assumptions of Theorem 29 hold true. If $f_x - x \in K^W$, then $\|Q_s^*\| = Q_s^*(f_x^-)$.*

*In addition, if $-f_x^- \in int(\mathcal{D})$ [or, equivalently, $-(1+\epsilon)f_x^- \in \mathcal{D}$ for some $\epsilon > 0$], then*

$$Q_s^* = 0$$

*and, consequently, $Q^* = Q_r^*$ is unique.*

PROOF. We know that $\|Q_s^*\| \geq Q_s^*(f_x^-)$. Replacing $k$ with $f_x - x \in K^W$, in (35)–(36), we obtain $E_{Q_r^*}[f_x] \leq Q_s^*(f_x^-) + x$. By the optimal relation (41), we necessarily have $\|Q_s^*\| = Q_s^*(f_x^-)$.

Suppose, now, that $-f_x^-$ also belongs to $int(\mathcal{D})$. Then, $-(1+\epsilon)f_x^- \in \mathcal{D}$ for some $\epsilon > 0$ and this implies that

$$Q_s^*(f_x^-) = \|Q_s^*\| \geq (1+\epsilon)Q_s^*(f_x^-),$$

whence $\|Q_s^*\| = 0$. □

This proposition provides a financial interpretation of the extra term $\|Q_s^*\|$ in (41): it is equal to the optimal singular part computed on the negative part of the optimal claim. Also, when the optimal loss $f_x^-$ is well inside the tolerated margin, the singular part is zero.

REMARK 32. In the case where $a$ is finite, $f_x - x$ is always in $K^1$, as we will see in Theorem 40, and so, in this case, $\|Q_s^*\| = Q_s^*(f_x^-)$.

Note that in the whole of Section 4, we have never required that the loss bound $W$ is also suitable (i.e., $W \in L^{\widehat{u}} \cap \mathbb{S}$). But, of course, the suitability condition, as given in Definition 2, is a desirable property. In fact, it guarantees that the domain of maximization $K^W$ is nontrivial.

Moreover, we have shown in (20) that when $W \in L^{\widehat{u}} \cap \mathbb{S}$, the regular measures in $\mathcal{M}^W$ can be characterized as the set $M_\sigma \cap L^{\widehat{\Phi}}$, *independently of* $W$. Since any dual optimum $Q^*$ has $Q^*_r \in \mathcal{L}^\Phi$, we deduce that, in the case $W \in L^{\widehat{u}} \cap \mathbb{S}$,

> when the singular part is null, $Q^*$ is unique and it is a $\sigma$-martingale measure with finite $\Phi$-entropy.

REMARK 33 (On no arbitrage). As already remarked in [5], Lemma 1, the hypothesis (38) in our main theorem is a totally different notion from No Free Lunch with Vanishing Risk, as defined by [13]. It does not ensure that there are $Q \in M_\sigma$ which are equivalent to $P$ and the existence of such equivalent measures does not imply (38).

A detailed analysis of the relationship between No Arbitrage or NFLVR and utility maximization (for a class of agents acting in the market) is provided in [16].

REMARK 34 (*Random endowments*). Our framework also allows for a complete treatment of the situation in which a random endowment is given, namely, one considers a maximization of the type

$$\sup_{k \in K^W} E[u(e+k)],$$

where $e$ is an $\mathcal{F}_T$ measurable r.v. satisfying some integrability properties. The paper [7] is entirely dedicated to the resolution of this problem.

## 5. Examples.

5.1. *Finite period markets with suitable and compatible loss bounds $W$.* In a finite period market, the filtration is formed by a finite number of increasing $\sigma$-algebras $(\mathcal{F}_t)_t$, $t = 0, 1, \ldots, T$. In this case, the analysis of the existence of suitable and weakly compatible loss bounds $W$ is rather straightforward.

Indeed, there is such a good $W \in L^{\widehat{u}} \cap \mathbb{S} \neq \varnothing$ (or $W$ suitable and compatible, i.e., $W \in M^{\widehat{u}} \cap \mathbb{S} \neq \varnothing$) iff for each $i = 1, \ldots, d$, there exists $H^i \in L(S^i)$ such that

$$H^i_t \neq 0 \quad \text{a.s. for all } t = 1, \ldots, T$$

and

$$|(H^i \cdot S^i)_t| \in L^{\widehat{u}} \text{ (or } M^{\widehat{u}}) \quad \text{for all } t = 1, \ldots, T.$$

Then, by setting $W \triangleq 1 + \sum_{i=1}^d (H^i \cdot S^i)^*_T$, we have $W \in L^{\widehat{u}}$ (resp. $M^{\widehat{u}}$) and it is obviously suitable [we denote with $Y^*$ the maximal process $(Y^*_t)_{t \geq 0}$, $Y^*_t = \sup_{s \leq t} |Y_s|$].

Clearly, if, for each $i = 1, \ldots, d$ and $t = 1, \ldots, T$, we have $|S^i_t| \in L^{\widehat{u}}$ (resp. $M^{\widehat{u}}$), then $W \triangleq 1 + \sum_{i=1}^d (S^i)^*_T \in L^{\widehat{u}} \cap \mathbb{S}$ (resp. $M^{\widehat{u}} \cap \mathbb{S}$).

5.2. *Exponential utility and $W \in (\widehat{L^u} \cap \mathbb{S}) \setminus \widehat{M^u}$.* We now give some examples of utility maximization problems that illustrate our results very well in the case $a = -\infty$ and in the novel situation $W \in (\widehat{L^u} \cap \mathbb{S}) \setminus \widehat{M^u}$. This is a case which was not covered by [9] and [5]. We consider the exponential utility $u(x) = -e^{-x}$ and assume that the initial endowment is zero in all of the examples.

EXAMPLE 35 $[-f_x^- \in int(\mathcal{D}), Q_s^* = 0]$. Let $S$ be a scalar compound Poisson process stopped at the finite horizon $T$, that is, $S_t = \sum_{(T_j \leq t \wedge T)} Y_j$, in which: (i) $T_0 = 0$ and $(T_j)_{j \geq 1}$ is the sequence of jump times of a Poisson process $N$ of parameter $\lambda$; (ii) $Y_0 = 0$ and $(Y_j)_{j \geq 1}$ is a sequence of i.i.d. random variables independent from $(T_j)_{j \geq 1}$ with doubly exponential distribution of parameter $\nu > 0$ and centered at 1 [i.e., the density is $f(y) = \frac{\nu}{2} e^{-\nu|y-1|}$]. This is the same kind of model we had in [5], the difference being that there, the $Y_j$ had a Gaussian distribution (for more details, see also [8]).

Then, $\mathcal{H}^1$ is trivial, while $W = 1 + \sup_{t \leq T} |S_t| \in \widehat{L^u} \cap \mathbb{S} \setminus \widehat{M^u}$. Maximizing over $\mathcal{H}^W$, we obtain

$$\sup_{k \in K^W} E[-e^{-k}] = \max_{k \in K_\Phi^W} E[-e^{-k}]$$

$$= \min_{\lambda > 0, Q \in \mathcal{N}^W} \left\{ \lambda E_{Q_r}\left[\ln\left(\frac{dQ_r}{dP}\right)\right] + \lambda(\ln \lambda - 1 + \|Q_s\|) \right\}.$$

The primal and dual optima are both unique:

$$f_0 = a^* S_T \in K^W,$$

$$Q^* = Q_r^*, \quad \frac{dQ^*}{dP} = \exp\left(-a^* S_T - \lambda T \left(\frac{\nu^2}{\nu^2 - (a^*)^2} e^{a^*} - 1\right)\right),$$

where $a^* = \sqrt{1 + \nu^2} - 1$. The uniqueness of $Q^* = Q_r^*$ follows from Proposition 31 since $f_0 \in K^W$ and $(1+\epsilon)f_0^- \in \mathcal{D}$ if $\epsilon < \frac{\nu - a^*}{a^*}$.

In the following three examples, we consider single period market models with $\mathcal{F}_0$ trivial. In this context, as soon as we exhibit a good loss bound $W$, then $\mathcal{H}^W = \mathbb{R}$.

EXAMPLE 36 $(Q_s^* \neq 0)$. Consider an exponential random variable $Y$ of parameter 1 [with density $f(y) = e^{-y} I_{[0,+\infty)}(y)$] and suppose that $\Omega$ supports a discrete r.v. $Z \in L^\infty$, independent from $Y$, which takes values in $\{1, -\frac{1}{2}, \ldots, \frac{1}{n} - 1, \ldots\}$. Then, define $S_0 = 0$ and $S_1 = ZY$. As shorthand, let $p_1 = P(Z = 1) > 0$ and $p_n = P(Z = \frac{1}{n} - 1) > 0, n \geq 2$.

The investor has exponential utility so that $W = 1 + Y \in L^{\widehat{u}}$, and this is clearly in $\mathbb{S}$. We then must consider

$$\sup_{h \in \mathbb{R}} E[-e^{-hS_1}].$$

It is not difficult to see that a necessary condition for the quantity to be maximized to be finite is that $-1 < h \leq 1$, which is also sufficient if we require that $p_n$ goes to zero very quickly. For now, suppose that the convergence speed of the $p_n$ is such that $E[S_1 e^{-hS_1}]$ is also finite for $-1 < h \leq 1$. Then, if $g(h) = E[-e^{-hS_1}]$, the derivative is $g'(h) = E[S_1 e^{-hS_1}]$ and, in case it is positive for all $-1 < h \leq 1$, the maximum is reached when $h = 1$. Note that

$$(42) \quad g'(h) = \sum_{n \geq 1} p_n z_n E[Y e^{-hz_n Y}] = p_1 E[Y e^{-hY}] + \sum_{n \geq 2} p_n z_n E[Y e^{-hz_n Y}]$$

so that

$$g'(h) \geq p_1 E[Y e^{-Y}] - \sum_{n \geq 2} p_n E[Y e^{-z_n Y}]$$

and the right term is strictly positive when the $(p_n)_{n \geq 2}$ are sufficiently small. So, we can assume that $g'(h) > 0$ for all $-1 < h \leq 1$.

In such a case, the optimal claim is $f_0 = S_1$ and the unique regular part of $Q^*$ is $\frac{dQ_r^*}{dP} = \frac{e^{-S_1}}{E[e^{-S_1}]}$. Since $g'(1) > 0$, we obtain

$$E_{Q_r^*}[S_1] = E\left[S_1 \frac{e^{-S_1}}{E[e^{-S_1}]}\right] = \frac{g'(1)}{E[e^{-S_1}]} > 0.$$

Hence, any optimal $Q^*$ necessarily has a nonzero singular part $Q_s^*$ since, by Proposition 31, we know that $E_{Q_r^*}[S_1] = Q_s^*(S_1^-) = \|Q_s^*\|$.

EXAMPLE 37. Here, we show that the condition $-f_x^- \in int(\mathcal{D})$ is only sufficient to obtain $Q_s^* = 0$. The setting is the same as that of the example above, up to equation (42). The $p_n$ can be selected so that $g'(h) > 0$ if $-1 < h < 1$, but $g'(1) = 0$. In fact, note that when $g'$ is positive, it is also monotone decreasing since $g''(h) = -p_1 E[Y^2 e^{-hY}] - \sum_{n \geq 1} p_n (z_n)^2 E[Y^2 e^{-hz_n Y}] < 0$. So, we only have to impose the condition $g'(1) = 0$, that is, $0 = p_1 E[Y e^{-Y}] + \sum_{n \geq 2} p_n z_n E[Y e^{-z_n Y}]$. In this way, $f_0$ is again $S_1$, but $E_{Q_r^*}[S_1] = 0$ and, by Proposition 31, $Q_s^* = 0$.

The next example is perhaps the most interesting since it is a concrete case in which the optimal functional $Q^*$ is not unique, there being an infinity of singular positive $Q_s \in \mathcal{A}^d$ such that $Q_r^* + Q_s \in \mathcal{M}^W$ and which satisfy the optimal relation (41).

EXAMPLE 38 (*An infinity of optimal functionals*). Let $(\Omega, \mathcal{F}, P)$ be the product space of two discrete spaces,

$$(\Omega_a = \{\omega_1^a, \omega_2^a, \ldots\}, \mathcal{F}_a, P_a) \quad \text{and} \quad (\Omega_b = \{\omega_1^b, \omega_2^b, \ldots\}, \mathcal{F}_b, P_b).$$

Also, suppose that $P_a(\omega_i^a) = P(\omega \mid \omega^a = \omega_i^a) > 0$ for all $i$ and that $P_b(\omega_j^b) = P(\omega \mid \omega^b = \omega_j^b) = (e-1)e^{-j}$. Now, let $W$ be the r.v. $W = \sum_{j \geq 1} j I_{\{\omega^b = \omega_j^b\}}$. Then, $P(W = j) = (e-1)e^{-j}$.

If $f \in L^0(\Omega)$, call $f_{ij} = f(\omega_i^a, \omega_j^b)$. Define the usual one-period market model on $\Omega$ ($\mathcal{F}_0$ trivial) as follows: $S_0 = 0$ and $S_1$ equal to

$$W I_{\{\omega^a = \omega_1^a\}} + \sum_{i > 1} -W I_{\{W \leq i\}} I_{\{\omega^a = \omega_i^a\}}$$

[think of $\Omega$ as a matrix $(\omega_i^a, \omega_j^b)_{ij}$: on the first row, $S_1$ is equal to $W$ and on the other rows, $S_1$ is equal to $-W$ up to the diagonal term $(\omega_i^a, \omega_i^b)$ and then null]. We can impose conditions on the $P_a(\omega_i^a) = P(\omega^a = \omega_i^a) > 0$ in order that $S_1 \in L^{\widehat{u}}$ and the resulting $f_x$ is equal to an arbitrary positive multiple of $S_1$ (say $5 S_1$).

Since we must again consider

$$\sup_{h \in \mathbb{R}} E[-e^{-h S_1}],$$

it is sufficient, as before, to show that we can require that:

1. $g(h) = E[-e^{-h S_1}]$ is finite iff $-1 < h \leq 5$;
2. $g'(h) > 0$ for $-1 < h \leq 5$.

We separately prove the two items above.

1. Note that $g(h) = p_1 E[-e^{-hW}] + \sum_{i>1} p_i E[-e^{hW I_{\{W \leq i\}}}]$ [where $p_i = P_a(\omega_i^a)$ for short] and that $h > -1$ is then obvious.

    To obtain $h \leq 5$, note that

    $$E[e^{h W I_{\{W \leq i\}}}] = (e-1) \sum_{k=1}^{i} e^{(h-1)k} + P(W > i),$$

    so that when $h > 1$ is fixed, this term for large $i$ is of the same order of magnitude as $e^{(h-1)i}$. If we select $p_1 = 1 - \sum_{i>1} p_i$ and $p_i \sim \frac{1}{i^r e^{4i}}$ with the power $r > 1$ arbitrary, we derive that $g(5)$ is finite, while $g(5 + \epsilon) = -\infty$.
2. Given these asymptotics, we show that, for some $r$,

    $$g'(h) = E[S_1 e^{-h S_1}]$$
    $$= p_1(r) E[W e^{-hW}] + \sum_{i>1} p_i(r) E[-W I_{\{W \leq i\}} e^{h W I_{\{W \leq i\}}}] > 0.$$

To this end, observe that the term

$$\sum_{i>1} p_i(r) E[W I_{\{W \leq i\}} e^{5WI_{\{W \leq i\}}}]$$

is infinitesimal when the power $r \to \infty$. In fact, $E[WI_{\{W \leq i\}} e^{5WI_{\{W \leq i\}}}] = (e-1)\sum_{k=1}^{i} k e^{5k} e^{-k} \leq (e-1) i^2 e^{4i}$. So, if $r > 3$,

$$\sum_{i>1} p_i(r) E[W I_{\{W \leq i\}} e^{5WI_{\{W \leq i\}}}] \leq C_1 \sum_{i>1} \frac{1}{i^r e^{4i}} i^2 e^{4i}$$

$$= C_1 \sum_{i>1} \frac{1}{i^{r-2}} = C_2 \frac{1}{r-3},$$

where $C_1, C_2$ are positive constants.

Hence, if $r$ is sufficiently large, $p_1 = p_1(r)$ is close to 1 and $g'(h) = E[S_1 e^{-hS_1}] = p_1 E[We^{-hW}] - \epsilon > 0$. We then have $\|Q_s^*\| = E_{Q_r^*}[f_x] = E[5 S_1 \frac{e^{-5S_1}}{E[e^{-5S_1}]}] > 0$ and $\|Q_s^*\| = Q_s^*(5 S_1^-)$.

Let us exhibit some different $Q_s^*$. To this end, we need the Hahn–Banach extension theorem. Consider the function $\psi$ defined on $\widehat{L^u}$ as

$$\psi(f) = \limsup_i \frac{f_{ii}}{i}.$$

It is then not difficult to show that $\psi$ is finite on $\widehat{L^u}$. In fact, $f \in \widehat{L^u}$ iff $E[e^{\alpha|f|}]$ is finite for some positive $\alpha$. Then,

$$\sum_i e^{\alpha|f_{ii}|} p_{ii} \leq \sum_{i,j} e^{\alpha|f_{ij}|} p_{ij} < +\infty$$

and $p_{ii} = P(\omega^a = \omega_i^a, \omega^b = \omega_i^b) = (e-1) e^{-i} p_i \sim \frac{1}{i^r e^{5i}}$. The convergence of the series $\sum_i e^{\alpha|f_{ii}|} \frac{1}{i^r e^{5i}}$ implies that the general term tends to 0, that is, $\lim_i \alpha|f_{ii}| - r\ln i - 5i = -\infty$. So, definitely,

$$\alpha |f_{ii}| < 5i + r \ln i \tag{43}$$

and $\psi(f)$ is finite. The function $\psi$ is evidently positively homogeneous, subadditive, null over $L^\infty$ and such that $\psi(-S_1) = 1$. Define $T_1$ to be the linear functional over $span(L^\infty, S_1)$ that is null on $L^\infty$ and such that $T_1(-S_1) = 1$. Since

$$T_1(f) = \psi(f),$$

by the Hahn–Banach theorem, $T_1$ can be extended to a linear functional $T \leq \psi$ on $\widehat{L^u}$.

In addition, $T$ is positive since, if $f \geq 0$, then $-T(f) = T(-f) \leq \psi(-f) \leq 0$. Namioka's theorem [2] ensures that $T$ is continuous. However, it is very easy to show continuity directly. We prove that $T$ is bounded on the positive elements of the open unit ball $B$ and henceforth continuous. In fact, $b \in B_+$ implies $E[\hat{u}(b)] \leq 1$ so that $E[e^b] < +\infty$ and, from (43), we know that definitely $b_{ii} < 5i + r \ln i$. Then,

$$\|T\| = \sup_{b \in B_+} T(b) \leq \sup_{b \in B_+} \psi(b) \leq 5$$

and it is almost immediate to show that $\|T\| = 5$. Since $T$ is a continuous functional which is null on $L^\infty$, it is null on $M^{\widehat{u}} = \overline{L^\infty}$. So, $T \in \mathcal{A}_+^d$ and if we positively scale it as

$$T_\psi(f) \triangleq E_{Q_r^*}[S_1] T(f),$$

then we have that $Q^* = Q_r^* + T_\psi$ is optimal since

$$Q^*(f_0) = Q^*(5 S_1) = 0.$$

The same argument can be repeated on the first subdiagonal of $\Omega$, $j = i - 1$ (or on any other subdiagonal), that is, one can consider

$$\varphi(f) = \limsup_i \frac{f_{i,i-1}}{i}$$

and construct the corresponding $T_\varphi$, which gives another $Q_s^*$, and so on.

## 6. Comparison with existing literature.

6.1. *The case $a$ finite.* In the case $a$ finite, one immediately thinks of the seminal paper [19]. In this article, the dual domain consists of terminal variables of nonnegative supermartingales $Y$ with $Y_0 = 1$. The authors pointed out that the dual optimum $Y_T^* \in L^1$ may not satisfy $E[Y_T^*] = 1$, but the approach in [19] does not provide an interpretation of the lost mass.

In [12], the authors are very concerned with this problem since they admit random endowments. In [12], the set of admissible strategies on which the maximization is performed is exactly $\mathcal{H}^1$ and so $C^1 = (K^1 - L_+^0) \cap L^\infty$. The dual variables in [12] are those in the set $(C^1)^0 \cap \{Q \in ba \mid Q(I_\Omega) = 1\} \subseteq ba$, which they call $\mathcal{D}$ and which coincides with $\mathcal{M}^1$, as defined in (17). Indeed, in our setting and when $a$ is finite, $L^{\widehat{u}} = L^\infty$, $(L^{\widehat{u}})^* = ba$ and the primal domain $K^1$ is obtained by selecting $W = 1 \in L^{\widehat{u}}$.

This is the reason why we prefer working out the more natural comparison with [12] instead of [19]. We show that if we set the random endowment $e = 0$ in [12], then it is rather easy to recover the known results thanks to our unifying Theorem 29. For a comparison of the results when $e \neq 0$, we refer to [7].

In [12] (as well as in [19]), $S$ is a possibly nonlocally bounded semimartingale, the utility function $u^0 : (0, +\infty) \to \mathbb{R}$ has $(0, +\infty)$ as proper domain and the assumptions are as follows:

(1) condition (A1) on $u^0$;
(2) $AE_{+\infty}(u^0) < 1$;
(3) there exists a probability $Q$ equivalent to $P$ such that for each $H \in \mathcal{H}^1$, the process $H \cdot S$ is a local martingale under $Q$ (NFLVR-type condition);
(4) $\sup_{k \in K^1} E[u^0(x+k)] < \infty$ for some $x > 0$.

(In [19], there is an irrelevant difference in the statement of (4), while in [12], this condition is equivalently formulated as $\sup_{k \in C^1} E[u^0(x+k)] < \infty$.) We now compare these assumptions with those adopted in this paper.

REMARK 39. Let us define
$$u(x) \triangleq u^0(x-a) \qquad \text{for some fixed } a < 0$$
so that the proper domain of $u$ is $(a, +\infty)$, as required in our paper. Then, the convex conjugate of $u^0$ is $\Phi^0(y) = \Phi(y) + ay$.
1. The hypothesis (A1) on $u^0$ clearly implies (A1) on $u$.
2. The hypothesis $AE_{+\infty}(u^0) < 1$ implies our assumptions (A2) on $\Phi$.
  Indeed, from Corollary 6.1(iii) [19], we know that the condition $AE_{+\infty}(u^0) < 1$ implies a "nice" behavior near zero of $\Phi^0$, and so also of $\Phi$. As shown in Section 2.1.1, in the case that $a$ is finite, we have $\Phi(y) \leq -ay + u(0)$ for large values of $y$ and therefore assumption (A2) on $\Phi$ holds true.
3. In Theorem 40 below, we will adopt the condition of NFLVR, which is equivalent to the assumption (3) above, as can be easily deduced from Theorem 1.1 and Proposition 4.7 [13]. Indeed, NFLVR is equivalent to the existence of a $Q \in \mathbb{M}_\sigma$ equivalent to $P$ and this implies that, for all $H \in \mathcal{H}^1$, the process $H \cdot S$ is a $Q$-local martingale [see Proposition 19(a), in case $W = 1$, for a proof of this well-known fact]. Thus, NFLVR implies the assumption (3) above. Conversely, note that $H \cdot S$ is uniformly bounded from below if $H \in \mathcal{H}^1$. Hence, (3) implies that there exists a probability $Q$ equivalent to $P$ such that, for all $H \in \mathcal{H}^1$, $H \cdot S$ is a $Q$-supermartingale, so that $E_Q[(H \cdot S)_T] \leq 0$. This also implies that $E_Q[f] \leq 0$ for all $f \in C$ and hence for all $f \in \overline{C}$, the $L^\infty$-closure of $C$. Hence, $\overline{C} \cap L^\infty_+ = \{0\}$ and NFLVR holds true.
4. When $a$ is finite and NFLVR holds true, the condition that appears in the next theorem—$\sup_{k \in K^1} E[u^0(x+k)] < \infty$ for some $x > 0$—is equivalent to the condition $\sup_{k \in K^1} E[u^0(x+k)] < u^0(\infty)$ for some $x > 0$ (see Remark 3.7 [24]) and therefore it is also equivalent to the condition that is used in this paper, $\sup_{k \in K^1} E[u(x+a+k)] < u(\infty)$ for some $x > 0$.

Then, from Theorem 29, we can derive the following, which sums up the results in [12] for the case with no random endowment.

THEOREM 40. *Suppose that NFLVR holds true. If $u^0$ satisfies* (A1), $AE_{+\infty}(u^0) < 1$ *and* $\sup_{k \in K^1} E[u^0(x+k)] < \infty$ *for some* $x > 0$, *then, for all* $x > 0$,

$$\text{(44)} \qquad \sup_{k \in K^1} E[u^0(x+k)] = \min_{\lambda > 0, Q \in \mathcal{M}^1} \lambda x + E\bigg[\Phi^0\bigg(\lambda \frac{dQ_r}{dP}\bigg)\bigg]$$

*and the solutions are related as*

$$\widehat{X} = -x - (\Phi^0)'\bigg(\lambda^* \frac{dQ_r^*}{dP}\bigg),$$

*where $\widehat{X} \in K^1$ is the primal solution (so that $\widehat{X}$ is replicable with a strategy $\widehat{H} \in \mathcal{H}^1$) and $\lambda^*$ (unique) and $Q^*$ (unique in the regular part) are the dual solutions. In addition, $\widehat{X}$ satisfies*

$$\text{(45)} \qquad E_{Q_r}[x + \widehat{X}] \le x \qquad \forall Q \in \mathcal{M}^1,$$

$$\text{(46)} \qquad E_{Q_r^*}[x + \widehat{X}] = x \qquad \text{for the optimal } Q^*.$$

PROOF. As explained in the above remark, we may apply Theorems 21 and 29 to the utility function $u$ defined by $u(x) \triangleq u^0(x - a)$ for some fixed $a < 0$.

From Theorem 21 and by (16),

$$\sup_{k \in K^1} E[u^0(x+k)] = \sup_{k \in K^1} E[u(x+a+k)]$$

$$= \min_{\lambda > 0, Q \in \mathcal{M}^1} \bigg\{\lambda(x+a) + E\bigg[\Phi\bigg(\lambda \frac{dQ_r}{dP}\bigg)\bigg] - a\lambda Q_s(\Omega)\bigg\}.$$

To obtain (44), simply substitute in the above relation $\Phi(\lambda \frac{dQ_r}{dP}) = \Phi^0(\lambda \frac{dQ_r}{dP}) - a\lambda \frac{dQ_r}{dP}$ and $Q(\Omega) = E[\frac{dQ_r}{dP}] + Q_s(\Omega) = 1$. From Theorem 29, we also know that the optimal $Q^*$ is unique in the regular part and that the link between primal and dual optima is

$$f_{(x+a)} = -\Phi'\bigg(\lambda^* \frac{dQ_r^*}{dP}\bigg) = -(\Phi^0)'\bigg(\lambda^* \frac{dQ_r^*}{dP}\bigg) + a.$$

Note that $f_{x+a} > a$. From (40), $f_{x+a} \in K^1_\Phi(x+a)$, whence

$$E_{Q_r}[f_{x+a}] \le x + a + \|Q_s\|$$

$$= (x+a) - aQ_s(\Omega) = x + aE_{Q_r}[I_\Omega] \qquad \forall Q \in \mathcal{N}^1,$$

$$E_{Q_r^*}[f_{x+a}] = x + a + \|Q_s^*\| = (x+a) - aQ_s^*(\Omega) = x + aE_{Q_r^*}[I_\Omega]$$

for the optimal $Q^*$. By definition, $f_{x+a} = (x+a) + \widehat{X}$. Therefore, the two relations above can be rewritten as

$$E_{Q_r}[x + \widehat{X}] \le x \qquad \forall Q \in \mathcal{N}^1,$$

$$E_{Q_r^*}[x + \widehat{X}] = x \qquad \text{for the optimal } Q^*,$$

from which we derive (46) and (45) for the set $\mathcal{N}^1$, but not yet for the larger set $\mathcal{M}^1 = \mathcal{D}$. But an easy approximation argument based on convexity, as in [12], Lemma 4.4, ensures that the above inequality holds for all $Q \in \mathcal{M}^1$. Finally, we show that $\widehat{X} \in K^1$. The NFLVR assumption implies, in particular, that $\mathcal{M}^1 \cap L^1 \neq \varnothing$. Since $f_{x+a} > a$, $\widehat{X}$ is also bounded from below, by $-x$, and $\widehat{X} \wedge n \in L^\infty$. From (45),

$$(47) \qquad E_Q[\widehat{X} \wedge n] \leq 0 \qquad \forall Q \in \mathcal{M}^1 \cap L^1.$$

Again by NFLVR, the convex cone $C^1$ is $\sigma(L^\infty, L^1)$-closed (see [13]). So, the relation (47), together with an application of the bipolar theorem, ensures that $\widehat{X} \wedge n$ (which is bounded from below by $-x$) belongs to $C^1$ for all $n$.

In [13], it was proven that the set

$$Z = \{g \in L^0 \mid \exists k \in K^1, k \geq -x \text{ and } g \leq k\}$$

is closed in probability. Since $Z \supseteq C^1$, $\widehat{X} \wedge n \in Z$ for all $n$ and

$$\widehat{X} \wedge n \uparrow_n \widehat{X},$$

we derive that $\widehat{X} \in Z$, so $\widehat{X} \leq \widehat{k}$ for some $\widehat{k} \in K^1$ and, by optimality, $\widehat{X} = \widehat{k} = (\widehat{H} \cdot S)_T$ a.s. $\square$

6.2. *The case $a = -\infty$ and $W \in M^{\widehat{u}}$.* The natural term of comparison is now [5] and we show that Theorem 29 permits the recovery of the results obtained there under the stronger hypothesis $W \in M^{\widehat{u}} \cap \mathbb{S}$.

Note that the condition $W \in M^{\widehat{u}} \cap \mathbb{S}$ already implies that $a = -\infty$ since, otherwise, $M^{\widehat{u}} = \{0\}$.

Recall that in the case $a = -\infty$, if $Q \in \mathcal{N}^W$, then $Q(I_\Omega) = E_{Q_r}[I_\Omega] = 1$ so that the regular parts of the elements in $\mathcal{N}^W$ are probability measures and, consequently,

$$K_\Phi^W(x) = x + K_\Phi^W(0).$$

LEMMA 41. *Let $W \in M_+^{\widehat{u}}$. If $z \in (C^W)^0$, then $z(t) \triangleq z_r + tz_s \in (C^W)^0$ for all $0 \leq t \leq 1$.*

PROOF. Each $f \in C^W$ can be written as $f = k \wedge g - h$, with $g, h \in L_+^{\widehat{u}}$ [select $k \in K^W$ so that $f \leq k$, then $f = k \wedge f^+ - (k \wedge f^+ - f)$]. Given this decomposition,

$$z_r(f) \leq z_r(k \wedge g) \leq z_r(k \wedge g) + z_s((k \wedge g)^+) = z(k \wedge g) \leq 0$$

$$\text{for all } f \in C^W, z \in (C^W)^0$$

since $W \in M^{\widehat{u}}$ implies $(k \wedge g)^- \in M^{\widehat{u}}$ (see Lemma 15) and, consequently, $z_s((k \wedge g)^-) = 0$. Therefore, $z_r \in (C^W)^0$. The result for $z(t)$ follows from the convexity of the polar. $\square$

If $W \in M_+^{\widehat{u}}$ and $Q \in \mathcal{N}^W$, then its regular part $Q_r$ is already in $\mathcal{N}^W$, that is,

$$\mathcal{N}^W \cap L^{\widehat{\Phi}} = \{Q_r \mid Q \in \mathcal{N}^W\}.$$

If, in addition, $W \in M^{\widehat{u}} \cap \mathbb{S}$, then, by (20) and the very definition of $\mathcal{N}^W$, the regular elements in $\mathcal{N}^W$ are $\mathbb{M}_\sigma \cap \mathbb{P}_\Phi$, independently of $W$. Therefore, we have the following corollary.

COROLLARY 42. *If $W \in M^{\widehat{u}} \cap \mathbb{S}$, then*

$$K_\Phi^W(0) = \left\{ f \in \bigcap_{Q \in \mathcal{N}^W} L^1(Q_r) \mid E_{Q_r}[f] \leq \|Q_s\| \text{ for all } Q \in \mathcal{N}^W \right\}$$

$$= \left\{ f \in \bigcap_{Q \in \mathbb{M}_\sigma \cap \mathbb{P}_\Phi} L^1(Q) \mid E_Q[f] \leq 0 \text{ for all } Q \in \mathbb{M}_\sigma \cap \mathbb{P}_\Phi \right\} \triangleq K_\Phi,$$

*where $K_\Phi$ is the domain used in [5], Theorem 1, and*

$$\min_{\lambda>0, Q \in \mathcal{N}^W} \lambda x + E\left[\Phi\left(\frac{dQ_r}{dP}\right)\right] + \lambda\|Q_s\|$$

$$= \min_{\lambda>0, Q \in \mathbb{M}_\sigma \cap \mathbb{P}_\Phi} \lambda x + E\left[\Phi\left(\frac{dQ}{dP}\right)\right] = U_\Phi(x).$$

Note that the domain of the primal problem ceases to depend on the particular $W$ selected as soon as $W \in M^{\widehat{u}} \cap \mathbb{S}$. Also, the dual problem reaches its minimal value on the set of probabilities $\mathbb{M}_\sigma \cap \mathbb{P}_\Phi$. Therefore, the dual can be reformulated so that no singular parts appear and the content of Theorem 29 coincides with the following result.

THEOREM 43 ([5], Theorem 1). *Suppose that assumptions* (A1) *and* (A2) *hold true and there exist $W \in M^{\widehat{u}} \cap \mathbb{S}$ and $x \in \mathbb{R}$ such that $U^W(x) < u(\infty)$. Then:*

(a) $\mathbb{M}_\sigma \cap \mathbb{P}_\Phi \neq \varnothing$;

(b) *for all $W \in M^{\widehat{u}} \cap \mathbb{S}$ and all $x \in \mathbb{R}$, the optimal value $U^W(x)$ is less than $u(\infty)$—it does not depend on the particular $W \in M^{\widehat{u}} \cap \mathbb{S}$ and*

$$U^W(x) = U_\Phi^W(x) = \min_{\lambda>0, Q \in \mathbb{M}_\sigma \cap \mathbb{P}_\Phi} \left\{ \lambda x + E\left[\Phi\left(\lambda \frac{dQ}{dP}\right)\right] \right\} = U_\Phi(x);$$

(c) *For all $x \in \mathbb{R}$, there exists the optimal solution $f_x = -x - \Phi'(\lambda_x \times \frac{dQ_x}{dP}) \in K_\Phi$,*

$$\max\{E[u(x+f)] \mid f \in K_\Phi\} = E[u(x+f_x)] = U_\Phi(x) < u(\infty),$$

*where $\lambda_x, Q_x$ are the optimal solution of the dual problem in item* (b).

**7. Which $W$?** Under the same assumptions as Theorem 43, we show, in the next proposition, that the optimal level of wealth that an investor may achieve by investing in $W$-admissible trading strategies, for *any* $W \in L_+^{\widehat{u}}$, is exactly $U_\Phi(x)$.

Of course, for a fixed $W \geq 1$ not necessarily suitable, $U^W(x)$ could be strictly less than $U_\Phi(x)$, as in Example 4 when $W = 1$.

We then derive that once an element $W_1 \in M^{\widehat{u}} \cap \mathbb{S}$ is identified, there is no incentive to invest in trading strategies $H \in \mathcal{H}^W$ with $W \in L^{\widehat{u}}$ and $W \geq W_1$.

Recall that $\mathcal{H}^{\widehat{u}} = \bigcup_{W \geq 1, W \in L^{\widehat{u}}} \mathcal{H}^W$.

PROPOSITION 44. *Suppose that assumptions* (A1) *and* (A2) *hold true and that there exist $W_1 \in M^{\widehat{u}} \cap \mathbb{S}$ and $x \in \mathbb{R}$ such that $U^{W_1}(x) < u(\infty)$. Then:*

1. *for all $W \in L_+^{\widehat{u}}$*

(48) $$U^W(x) \leq U_\Phi(x);$$

2.

(49)
$$\sup_{H \in \mathcal{H}^{\widehat{u}}} E[u(x + (H \cdot S)_T)]$$
$$= \min_{\lambda > 0, Q \in \mathbb{M}_\sigma \cap \mathbb{P}_\Phi} \left\{ \lambda x + E\left[\Phi\left(\lambda \frac{dQ}{dP}\right)\right] \right\} = U_\Phi(x);$$

3. *if $W \in L^{\widehat{u}}$ is greater than some $\overline{W} \in M^{\widehat{u}} \cap \mathbb{S}$, then*

$$U^W(x) = U^{\overline{W}}(x) = U_\Phi(x)$$

*and there is no incentive to invest in the strategies in $\mathcal{H}^W$.*

PROOF. From Remark 23 and Proposition 19(b), we have

$$\mathbb{M}_\sigma \cap \mathbb{P}_\Phi \subseteq \{Q \ll P \mid Q \in L^{\widehat{\Phi}} \text{ and } H \cdot S \text{ is a } Q\text{-supermartingale } \forall H \in \mathcal{H}^{\widehat{u}}\}.$$

Let $W \in L_+^{\widehat{u}}$. We then deduce (48) from the inequalities

$$U^W(x) = \sup_{k \in K^W} E[u(x+k)] \leq \sup_{k \in K^{(W+1)}} E[u(x+k)]$$

$$\leq \sup_{H \in \mathcal{H}^{\widehat{u}}} E[u(x + (H \cdot S)_T)]$$

$$\leq \inf_{\lambda > 0, Q \in \mathbb{M}_\sigma \cap \mathbb{P}_\Phi} \left\{ \lambda x + E\left[\Phi\left(\lambda \frac{dQ}{dP}\right)\right] \right\} = U_\Phi(x),$$

where the last inequality comes from the Fenchel inequality

$$u(x + (H \cdot S)_T) \leq \lambda \frac{dQ}{dP}(x + (H \cdot S)_T) + \Phi\left(\lambda \frac{dQ}{dP}\right)$$

and from $E_Q[(H \cdot S)_T] \leq 0$.

To show (49), let $W_1 \in M^{\widehat{u}} \cap \mathbb{S}$. We may apply Theorem 43 so that $\mathbb{M}_\sigma \cap \mathbb{P}_\Phi$ is not empty and, trivially,

$$\sup_{H \in \mathcal{H}^{\widehat{u}}} E[u(x + (H \cdot S)_T)] \geq \sup_{k \in K^{W_1}} E[u(x + k)] = U_\Phi(x).$$

Equality must then hold due to the opposite inequality given by (48).

Finally, if $W \in L^{\widehat{u}}$ and $W \geq \overline{W}$ for some $\overline{W} \in M^{\widehat{u}} \cap \mathbb{S}$, then $K^W \supseteq K^{\overline{W}}$ and

$$U^W(x) = \sup_{k \in K^W} E[u(x+k)] \geq \sup_{k \in K^{\overline{W}}} E[u(x+k)] = U_\Phi(x)$$

so that $U^W(x) = U^{\overline{W}}(x)$.  □

REMARK 45. When $M^{\widehat{u}} \cap \mathbb{S} \neq \varnothing$, we may directly state the primal optimization problem over the domain $(K^W - L_+^0) \cap M^{\widehat{u}}$ (see Lemma 15). The dual variables then live in the space $(M^{\widehat{u}})^* = L^{\widehat{\Phi}}$ so that no singular component appears and the results in Theorem 43 can be recovered by applying the duality between $M^{\widehat{u}}$ and $L^{\widehat{\Phi}}$. This is exactly the approach adopted in [9].

## APPENDIX

The representation of the conjugate of a convex integral functional on Orlicz spaces is provided by [18], Theorem 2.6 and is based on the similar representation on the space $L^\infty$, proven by [22]. In our notation, this theorem can be restated as follows.

THEOREM 46. *Suppose that $F : \mathbb{R} \to (-\infty, +\infty]$ and $F^* : \mathbb{R} \to (-\infty, +\infty]$ are convex l.s.c. functions (not identically equal to $+\infty$) conjugate to each other and that there exists $f \in L^{\widehat{u}}$ such that $I_F(f) \triangleq E[F(f)] < \infty$. If $I_{F^*}(g) <$* $\infty$ for some $g \in L^{\widehat{\Phi}}$, then the convex conjugate $I_F^*\colon (L^{\widehat{u}})^* \to (-\infty, +\infty]$ of the convex integral functional $I_F$ is given by

$$I_F^*(z) = I_{F^*}\left(\frac{dz_r}{dP}\right) + \sup\{z_s(f) \mid f \in \mathrm{dom}(I_F)\},$$

where $\mathrm{dom}(I_F)$ is the proper domain of $I_F$.

## REFERENCES

[1] ALIPRANTIS, C. and BURKINSHAW, O. (1985). *Positive Operators*. Academic Press, Orlando, FL. MR0809372
[2] ALIPRANTIS, C. and BORDER, K. C. (2005). *Infinite Dimensional Analysis*, 3rd ed. Springer, Berlin. MR1321140
[3] ANDO, T. (1960). Linear functionals on Orlicz spaces. *Nieuw Arch. Wisk. (3)* **8** 1–16. MR0123907
[4] ANSEL, J. P. and STRICKER, C. (1994). Couverture des actifs contingents et prix maximum. *Ann. Inst. H. Poincaré Probab. Statist.* **30** 303–315. MR1277002
[5] BIAGINI, S. and FRITTELLI, M. (2005). Utility maximization in incomplete markets for unbounded processes. *Finance and Stochastics* **9** 493–517. MR2212892
[6] BIAGINI, S. and FRITTELLI, M. (2007). The supermartingale property of the optimal portfolio process for general semimartingales. *Finance and Stochastics* **11** 253–266. MR2295831
[7] BIAGINI, S., FRITTELLI, M. and GRASSELI, M. (2007). Indifference price for general semimartingales. Unpublished manuscript.
[8] BIAGINI, S. (2005). Convex duality in financial theory with general semimartingales. Ph.D. thesis, Scuola Normale Superiore, Pisa.
[9] BIAGINI, S. (2008). An Orlicz spaces duality for utility maximization in incomplete markets. In *Seminar on Stochastic Analysis, Random Fields and Applications V. Progress Probab.* **59** 445–455. Birkhäuser, Basel.
[10] BISMUT, J. M. (1973). Conjugate convex functions in optimal stochastic control. *J. Math. Anal. Appl.* **44** 384–404. MR0329726
[11] BREZIS, H. (1983). *Analyse Fonctionnelle*. Masson, Paris. MR0697382
[12] CVITANIĆ, J., SCHACHERMAYER, W. and WANG, H. (2001). Utility maximization in incomplete markets with random endowment. *Finance and Stochastics* **5** 237–259. MR1841719
[13] DELBAEN, F. and SCHACHERMAYER, W. (1998). The fundamental theorem of asset pricing for unbounded stochastic processes. *Math. Ann.* **312** 215–250. MR1671792
[14] EMERY, M. (1980). Compensation de processus à variation finie non localement intégrables. *Séminaire de Probabilités XIV (Paris, 1978/1979). Lecture Notes in Math.* **784** 152–160. Springer, Berlin. MR0580120
[15] EKELAND, I. and TEMAM, R. (1976). *Convex Analysis and Variational Problems*. North-Holland, Amsterdam. MR0463994
[16] FRITTELLI, M. (2004). Some remarks on arbitrage and preferences in securities market models. *Math. Finance* **14** 351–357. MR2070168
[17] KAKUTANI, S. (1941). Concrete representation of abstract (L)-spaces and the mean ergodic theorem. *Ann. of Math.* **42** 523–537. MR0004095
[18] KOZEK, A. (1979). Convex integral functionals on Orlicz spaces. *Ann. Soc. Math. Polonae Ser. 1. Comm. Math.* **XXI** 109–134. MR0577677

[19] KRAMKOV, D. and SCHACHERMAYER, W. (1999). The asymptotic elasticity of utility function and optimal investment in incomplete markets. *Ann. Appl. Probab.* **9** 904–950. MR1722287

[20] PLISKA, S. R. (1986). A stochastic calculus model of continuous trading: Optimal portfolio. *Math. Oper. Res.* **11** 371–382. MR0844013

[21] RAO, M. M. (1968). Linear functionals on Orlicz spaces: General theory. *Pacific J. Math.* **25** 553–585. MR0412791

[22] ROCKAFELLAR, R. T. (1971). Integrals which are convex functionals. II. *Pacific J. Math.* **39** 694–734. MR0310612

[23] RAO, M. M. and REN, Z. D. (1991). *Theory of Orlicz Spaces*. Dekker, New York. MR1113700

[24] SCHACHERMAYER, W. (2001). Optimal investment in incomplete markets when wealth may become negative. *Ann. Appl. Probab.* **11** 694–734. MR1865021

DEPARTMENT OF ECONOMICS,
FINANCE AND STATISTICS
UNIVERSITY OF PERUGIA
VIA PASCOLI 20
PERUGIA 06123
ITALY
E-MAIL: s.biagini@unipg.it

DEPARTMENT OF MATHEMATICS
UNIVERSITY OF MILAN
VIA CESARI SALDINI 50
MILAN 20133
ITALY
E-MAIL: marco.frittelli@mat.unimi.it